\let\rup=\rightharpoonup % symbol for nontrivial action

\title{Regular embeddings of complete bipartite graphs: classification and enumeration}

\author{Gareth A. Jones\\
School of Mathematics\\
University of Southampton\\
Southampton SO17  1BJ, U.K.\\
{\tt G.A.Jones@maths.soton.ac.uk}
}

\documentclass[12pt]{article}
\usepackage[latin1]{inputenc}
\usepackage{a4wide}
\usepackage{latexsym}
\usepackage{amsmath}
\usepackage{amssymb}
\usepackage{amsfonts}
\newtheorem{thm}{Theorem}[section]
\newtheorem{lemma}[thm]{Lemma}
\newtheorem{cor}[thm]{Corollary}
\newtheorem{prop}[thm]{Proposition}
\date{}
\begin{document} 

\maketitle

\begin{abstract}
\noindent The regular embeddings of complete bipartite graphs $K_{n,n}$ in orientable surfaces are classified and enumerated, and their automorphism groups and combinatorial properties are determined. The method depends on earlier classifications in the cases where $n$ is a prime power, obtained in collaboration with Du, Kwak, Nedela and \v Skoviera, together with results of It\^o, Hall,
Huppert and Wielandt on factorisable groups and on finite solvable groups.
\end{abstract}



{\bf Running head:} Regular embeddings\\

\section{Introduction}

A major problem in topological graph theory is that of classifying the regular embeddings, in
orientable or non-orientable surfaces, of particular classes of arc-transitive graphs. Extra motivation has been provided by Grothendieck's theory of {\it dessins d'enfants\/}~\cite{Gro}, in which bipartite maps on orientable surfaces correspond to algebraic curves defined over the field $\overline{\bf Q}$ of algebraic numbers. In recent years there has been particular interest in the regular embeddings of the complete bipartite graphs $K_{n,n}$. The non-orientable regular embeddings of these graphs have recently been classified by Kwak and Kwon~\cite{KK3}, and our aim here is to present a similar result for the orientable regular embeddings. By these we mean the embeddings $\cal M$ of $K_{n,n}$ in orientable surfaces such that the orientation-preserving automorphism group ${\rm Aut}^+{\cal M}$ acts transitively on the arcs
(directed edges) of $\cal M$; for conciseness, we will refer to these simply as regular embeddings.

\begin{picture}(100,120)(-270,-40)

\put(-250,18){(a)}

\put(-140,20){\circle*{5}}
\put(-200,20){\circle{5}}
\put(-183,45){\circle*{5}}
\put(-155,45){\circle{5}}
\put(-183,-5){\circle*{5}}
\put(-155,-5){\circle{5}}

\put (-183,45){\line(1,0){26}}
\put (-183,-5){\line(1,0){26}}
\put (-140,20){\line(1,0){14}}
\put (-215,20){\line(1,0){13}}
\put (-199,21){\line(2,3){15}}
\put (-154,-4){\line(2,3){15}}
\put (-154,47){\line(2,3){8}}
\put (-192,-18){\line(2,3){7}}
\put (-199,18){\line(2,-3){15}}
\put (-155,43){\line(2,-3){15}}
\put (-154,-7){\line(2,-3){7}}
\put (-192,59){\line(2,-3){8}}

\multiput(-216,-5)(0,3){17}{.}
\multiput(-127,-5)(0,3){18}{.}
\multiput(-216,45)(2.55,1.5){18}{.}
\multiput(-170,-30)(2.55,1.5){17}{.}
\multiput(-216,-5)(2.55,-1.5){18}{.}
\multiput(-170,72)(2.55,-1.5){17}{.}

\put(-40,18){(b)}

\put(40,60){\circle*{5}}
\put(20,40){\circle*{5}}
\put(40,40){\circle{5}}
\put(60,40){\circle*{5}}
\put(0,20){\circle*{5}}
\put(20,20){\circle{5}}
\put(40,20){\circle*{5}}
\put(60,20){\circle{5}}
\put(80,20){\circle*{5}}
\put(20,0){\circle*{5}}
\put(40,0){\circle{5}}
\put(60,0){\circle*{5}}
\put(40,-20){\circle*{5}}

\put (1,20){\line(1,0){18}}
\put (21,20){\line(1,0){18}}
\put (41,20){\line(1,0){18}}
\put (61,20){\line(1,0){18}}
\put (21,40){\line(1,0){18}}
\put (41,40){\line(1,0){18}}
\put (21,0){\line(1,0){18}}
\put (41,0){\line(1,0){18}}
\put (20,21){\line(0,1){18}}
\put (40,21){\line(0,1){18}}
\put (60,21){\line(0,1){18}}
\put (20,1){\line(0,1){18}}
\put (40,1){\line(0,1){18}}
\put (60,1){\line(0,1){18}}
\put (40,-19){\line(0,1){18}}
\put (40,41){\line(0,1){18}}
\multiput(0,20)(2,2){20}{.}
\multiput(40,-20)(2,2){20}{.}
\multiput(-1,20)(2,-2){20}{.}
\multiput(39,60)(2,-2){20}{.}

\end{picture}

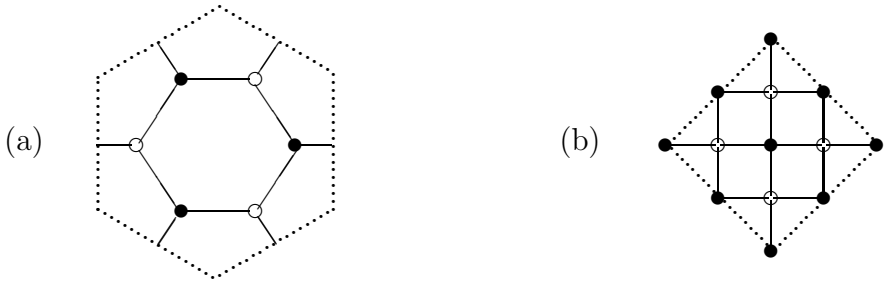
\centerline{Figure 1: embeddings ${\cal S}(3)$ and ${\cal N}(4\,;0,0)$; opposite sides of}
\centerline{the outer hexagon and square are identified to form a torus.}

\medskip

There is at least one such embedding of $K_{n,n}$ for each positive integer $n$, namely the standard embedding ${\cal S}(n)$, introduced by Biggs and White~\cite[\S 5.6.7]{BW} as a Cayley map for the group ${\bf Z}_{2n}$ with respect to the generators $1, 3, \ldots, 2n-1$ in that cyclic order. (By contrast, Biggs~\cite{Big} showed that the complete graph $K_n$ has regular embeddings, subsequently classified by James and the present author in~\cite{JJ}, if and only if $n$ is a prime power.) This embedding ${\cal S}(n)$, illustrated for $n=3$ in Fig.~1(a),  has $2n$-gonal faces, and has genus $(n-1)(n-2)/2$; it also appears in the theory of {\it dessins d'enfants\/} as a map on the $n$th degree Fermat curve $x^n+y^n=z^n$, with the black and white vertices and the edges represented by the inverse images of the points $0$, $1$ and the interval $[0,1]$ under the function $[x,y,z]\mapsto (x/z)^n$; see~\cite{Jon, JS} for details, and~\cite{CJSW, JSW} for some generalisations.

The classification process was begun by Nedela, \v Skoviera and Zlato\v s~\cite{NSZ}, who showed that if $n$ is prime then the standard embedding is the only regular embedding of $K_{n,n}$; more generally, it was shown in~\cite{JNS1}, using group theory, that $n$ has this uniqueness property if and only if $n$ is coprime to $\phi(n)$, where $\phi$ is Euler's function. (It follows from results of Burnside~\cite{Bur} and H\"older~\cite{Hol} that these are also the integers $n$ for which the cyclic group is the only group of order $n$~\cite[Exercise~575]{Ros}; Erd\H os~\cite{Erd} has shown that the proportion of the natural numbers $n\leq N$ with this property behaves asymptotically like
\[\frac{e^{-\gamma}}{\log\log\log N}\]
as $N\to\infty$, where $\gamma$ is Euler's constant.) Using
combinatorial methods based on permutations, Kwak and Kwon~\cite{KK1} have classified the regular embeddings in the case where $n$ is a product of two primes, and also the reflexible embeddings (those regular embeddings with an orientation-reversing automorphism) for all $n$~\cite{KK2}. The classification problem has also been solved for prime powers $n=p^e$: the regular embeddings for $p>2$ are classified in~\cite{JNS2} using properties of certain metacyclic $p$-groups, their direct analogues for $p=2$, associated with metacyclic $2$-groups, are dealt with in~\cite{DJKNS1}, and~\cite{DJKNS2} describes the small number of exceptional embeddings ${\cal N}(n\,;k,l)$, defined in \S 3, which arise for each $n=2^e\geq 4$ (see Fig.~1(b) for ${\cal N}(4\,;0,0)$). These prime power classifications, which form an important ingredient in the present paper, are summarised in \S 3, and they have been used by Coste, Streit, Wolfart and the present author in~\cite{CJSW, JSW} to study the action of the absolute Galois group ${\rm Gal}\,\overline{\bf Q}/{\bf Q}$ on the corresponding algebraic curves.

Our aim here is to show that for each $n$, the regular embeddings of $K_{n,n}$ can be
constructed from regular embeddings of $K_{p^e,p^e}$ for various prime powers $p^e$, using two simple
constructions described in \S4. The first of these, introduced in~\cite{JNS1}, is the cartesian
product, which produces a regular embedding $\cal S$ of $K_{s,s}$ from regular embeddings ${\cal
S}_i$ of $K_{s_i,s_i}\;(i=1,\ldots, k)$ for mutually coprime integers $s_i$ with $s=s_1\ldots s_k$. The
second construction produces a regular embedding $\cal M$ of $K_{n,n}$ as a $t^2$-sheeted regular
abelian covering of a regular embedding $\cal S$ of $K_{s,s}$, where $n=st$ with $s$ and $t$ coprime.
Our main result, formulated more precisely as Corollary~4.2 in \S 4, states that every regular
embedding of $K_{n,n}$ can be formed in this way as an abelian covering, where $\cal S$ is the
cartesian product of regular embeddings of $K_{s_i,s_i}$ for various prime powers $s_i$ appearing in
the factorisation of $n$, and $t$ is the product of the remaining prime powers. Since the prime
power regular embeddings are all known, this completes the classification of regular embeddings of
$K_{n,n}$ for all $n$. The proof is given in \S 8, after some preliminary results in \S\S5--7. We consider
the type and genus of these maps in \S9, and their mirror images and Petrie duals in \S10, while
formul\ae\/ for enumerating them are developed in \S11 and \S12. There are comments on connections with directed graphs and the Erd\H os-R\'enyi random graph in \S 13 and \S 14.

Our method is almost entirely group-theoretic, based on a study of the group ${\rm Aut}^+_0{\cal M}$
of automorphisms of a regular embedding $\cal M$ of $K_{n,n}$, preserving surface orientation and
vertex colours. We use the fact, proved in~\cite{JNS1} and briefly explained in \S 2, that a group
$G$ is isomorphic to ${\rm Aut}^+_0{\cal M}$ for some regular embedding $\cal M$ of $K_{n,n}$ if and
only if $G$ is isobicyclic, that is, a product of two disjoint cyclic groups of order $n$, with
an automorphism of $G$ transposing their generators. A result of Wielandt~\cite{Wie} on products of
nilpotent groups shows that such a group $G$ has a series
\[1=N_0<N_1<\cdots<N_{l-1}<N_l=G\]
of characteristic subgroups $N_i$ in which the quotients $N_i/N_{i-1}$ are isomorphic to the Sylow
$p$-subgroups $P$ of $G$ for the distincts primes $p$ dividing $n$. It follows that $\cal M$ can be
formed from a sequence of regular coverings, with the Sylow subgroups as the covering groups. These
Sylow subgroups $P$ are also isobicyclic, so they correspond to regular embeddings $\cal P$ of
$K_{p^e,p^e}$ where $p^e$ ranges over the prime powers in the factorisation of $n$. The classification
problem for prime powers having been solved, it is sufficient to consider how the groups $P$ and
maps $\cal P$ corresponding to the prime powers dividing $n$ may be combined to form $G$ and $\cal
M$. In particular, we need to determine how each Sylow $q$-subgroup $Q$ of $G$ can act by conjugation
on the subgroups $N_i$ to induce automorphisms of a Sylow $p$-subgroup $P\cong N_i/N_{i-1}$. Each
possible pattern of actions can be represented as a directed graph $\Gamma$, in which the vertices are
the primes $p$ dividing $n$, and an arc from $q$ to $p$ represents a nontrivial action of $Q$
on $P$. By determining the possible structures of such graphs $\Gamma$, and the ways in which regular
embeddings $\cal P$ and actions $Q\to{\rm Aut}\,P$ can be associated with their vertices and arcs, we
obtain a classification of isobicyclic groups $G$; this is stated as Theorem~4.1 in \S 4, and proved in
\S8. In particular, we find that the Sylow subgroups $P$ corresponding to primes $p$ which are terminal
vertices of arcs of $\Gamma$ generate a normal subgroup $T\cong C_t\times C_t$ of $G$, while for the
remaining primes dividing $n$ one can choose Sylow subgroups $S_i$ to generate their direct product $S$, which is a complement for $T$ in $G$; the action of $S$ by conjugation on $T$ can also be described explicitly, so that the structure of $G$ is completely known. Corollary~4.2, which classifies the regular embeddings of $K_{n,n}$, is an immediate consequence of Theorem~4.1, using the correspondence between groups and maps explained in \S2: the complement $S=S_1\times\cdots\times S_k$ corresponds to the cartesian product ${\cal S}={\cal S}_1\times\cdots\times{\cal S}_k$ in the classification theorem, with the normal subgroup $T$ yielding the abelian covering.

This method allows us to classify the reflexible embeddings of $K_{n,n}$, and those which are
self-Petrie (isomorphic to their Petrie duals), thus confirming the enumerations of these two classes
of maps recently obtained by Kwak and Kwon~\cite{KK2}. More generally, the method also yields an enumeration of all the regular embeddings of $K_{n,n}$. The general formula, which involves summation over all directed graphs $\Gamma$ which may occur for a given $n$, is rather complicated, and this is particularly so when $n$ is even because of the exceptional embeddings of $K_{2^e,2^e}$; however, in special cases, such as when $n$ is divisible by at most two primes, the formula simplifies to something relatively straightforward. In particular, the formula agrees with the enumerations obtained earlier in the special cases mentioned above.

\section{Regular embeddings and isobicyclic groups}

There is a detailed treatment in~\cite{JNS1} of the group-theoretic approach to regular embeddings of
complete bipartite graphs, so we will simply outline it here. Throughout this paper, $\cal M$ will
denote a regular embedding of $K_{n,n}$ in an orientable surface, and ${\rm Aut}^+{\cal M}$ will
denote its orientation-preserving automorphism group. The vertices can be coloured black or white so
that each edge connects vertices of different colours; then ${\rm Aut}_0^+{\cal M}$ denotes the
subgroup of index $2$ in ${\rm Aut}^+{\cal M}$ preserving the vertex colours. It is shown in~\cite{JNS1}
that a group $G$ is isomorphic to ${\rm Aut}^+_0{\cal M}$ for some such $\cal M$ if and only if $G$ is
$n$-{\sl isobicyclic\/}, or simply {\sl isobicyclic\/}; this means that $G$ has cyclic subgroups
$X=\langle x\rangle$ and $Y=\langle y\rangle$ of order $n$ such that $G=XY$ and $X\cap Y=1$ (so
$|G|=n^2$), and there is an automorphism $\alpha$ of $G$ transposing $x$ and $y$. We call $(G,x,y)$ an
$n$-{\sl isobicyclic triple\/}, and we call $x$ and $y$ the {\sl canonical generators\/} of $G$: these
are the orientation-preserving automorphisms of $\cal M$ fixing a black vertex $v$ and a white vertex
$w$, and sending each neighbouring vertex to the next vertex by following the orientation around $v$ or
$w$; the automorphism $\alpha$ of $G$ is induced by conjugation by the element of ${\rm Aut}^+{\cal M}$
which reverses the edge $vw$, so that ${\rm Aut}^+{\cal M}$ is a semidirect product
$G:\negthinspace\langle\alpha\rangle$ of $G$ by $C_2$. Conversely, given an $n$-isobicyclic triple $(G,
x, y)$, we can take the black and white vertices to be the cosets $gX$ and $gY$ of $X=\langle x\rangle$
and $Y=\langle y\rangle$ in $G$, and the edges to be the elements of $G$, with incidence given by
inclusion; the successive powers of $x$ and $y$ give the rotation of edges around each vertex, thus
defining the faces of a regular embedding $\cal M$ of $K_{n,n}$.

We define two isobicyclic triples $(G, x, y)$ and $(G', x', y')$ to be {\sl isomorphic\/} if there is
an isomorphism $G\to G'$ sending $x$ to $x'$ and $y$ to $y'$. Two regular embeddings of $K_{n,n}$ are
isomorphic if and only if their corresponding isobicyclic triples are isomorphic, so we have the following result:

%Theorem 2.1.

\begin{thm} The mapping ${\cal M}\mapsto(G,x,y)$ induces a bijection between the set
${\cal R}(n)$ of isomorphism classes of regular embeddings of $K_{n,n}$ and the set ${\cal
I}(n)$ of isomorphism classes of $n$-isobicyclic triples. Here ${\rm Aut}_0^+{\cal M}\cong G$ and
${\rm Aut}^+{\cal M}\cong G:\negthinspace\langle\alpha\rangle$. \hfill$\square$
\end{thm}

We will denote the standard embedding of $K_{n,n}$ by ${\cal S}(n)$; this corresponds to the {\sl
standard triple\/} $(G,x,y)\in{\cal I}(n)$, with $G=\langle x, y\mid x^n=y^n=[x,y]=1\rangle\cong
C_n\times C_n$. For any property of groups (such as being abelian), we will say that a map ${\cal
M}\in{\cal R}(n)$ and its corresponding triple $(G, x, y)\in{\cal I}(n)$ have this property if and
only if $G$ has this property. For instance, it is easily seen that ${\cal S}(n)$ is the only abelian
map in ${\cal R}(n)$.

\section{The prime power embeddings}

Here we briefly summarise the classifications in~\cite{JNS2, DJKNS1} and~\cite{DJKNS2} of what we
will call the {\sl prime power\/} embeddings and triples, by which we mean the regular embeddings
$\cal M$ of $K_{n,n}$ for prime powers $n$ and the assocated isobicyclic triples $(G, x, y)$. For our
purposes, it is sufficient to describe the groups $G={\rm Aut}_0^+{\cal M}$ and to give
representatives of the orbits of ${\rm Aut}\,G$ on pairs of canonical generators $x, y$.

As shown in~\cite{JNS2}, if $n=p^e$ for an odd prime $p$ then $G$ is a metacyclic group
$$G_f=\langle g,h\mid g^n=h^n=1,\,h^g=h^{1+p^f}\rangle \eqno(3.1)$$
where $f=1,\ldots, e$. This is a semidirect product $C_n\negthinspace :\negthinspace C_n$ of a normal
subgroup $\langle h\rangle\cong C_n$ by a complement $\langle g\rangle\cong C_n$, and the action of
$\langle g\rangle$ by conjugation on $\langle h\rangle$ is determined by the parameter $f$, which we
will sometimes denote by $f_G$. The canonical generators of $G$ can be chosen to be $x=g^u$ and
$y=g^uh$ for some $u$ coprime to $p$; we will denote the corresponding regular embedding and
isobicyclic triple by ${\cal M}(n,f,u)$ and ${\cal I}(n,f,u)$. Then ${\cal M}(n,f,u)\cong{\cal
M}(n,f',u')$, or equivalently ${\cal I}(n,f,u)\cong{\cal I}(n,f',u')$, if and only if $f=f'$ and
$u\equiv u'$ mod~$(p^{e-f})$, so there are $\phi(p^{e-f})$ isomorphism classes of regular embeddings
${\cal M}(n,f,u)$ for each $f$, represented by taking $u=1,\ldots, p^{e-f}$ coprime to $p$, and these
give a total of $\sum_{f=1}^e\phi(p^{e-f})=p^{e-1}$ maps in ${\cal R}(n)$ for each odd prime power
$n=p^e$. They all have type $\{2n,n\}$ and genus $(n-1)(n-2)/2$. Note that $G_e\cong C_n\times C_n$, and the unique map corresponding to $G_e$ is the standard
embedding ${\cal S}(n)={\cal M}(n,e,1)$.

It is shown in~\cite{DJKNS1} that if $n=2^e$ then the metacyclic maps ${\cal M}\in{\cal R}(n)$ and
triples $(G, x, y)\in{\cal I}(n)$ are the direct analogues for $p=2$ of those defined in $(3.1)$, but
with the case $f=1$ excluded if $e\geq 2$; this gives $\sum_{f=2}^e\phi(2^{e-f})=2^{e-2}$ metacyclic
maps and triples if $n=2^e\geq 4$, and one if $n=2$. In addition, it is shown in~\cite{DJKNS2} that there
is one non-metacyclic triple in ${\cal I}(4)$, and there are four in ${\cal I}(n)$ for each $n=2^e\geq
8$. To construct these, let $G(n\,;k,l)$ be the group defined by the presentation
\begin{eqnarray*}
G=G(n;k,l)=\langle x, y\;\mid & x^n=y^n=1,\, c:=[y,x]=x^{2+k2^{e-1}}y^{-2-k2^{e-1}},\\
& c^x=c^{-1+l2^{e-2}}x^4,\, c^y=c^{-1-l2^{e-2}}y^{-4}\rangle
\end{eqnarray*}
% \eqno(3.2)
where $k,l\in\{0, 1\}$. This is an extension of a normal subgroup $\Phi(G)=\langle x^2, y^2\rangle\cong
C_{n/2}\times C_{n/2}$ (the Frattini subgroup of $G$) by $G/\Phi(G)\cong C_2\times C_2$; the
conjugation action is given by $(y^2)^x=y^{-2}z^l$ and $(x^2)^y=x^{-2}z^l$ where $z$ is the central
involution $x^{n/2}y^{n/2}$. For each $n=2^e\geq 8$ the four choices for $k$ and $l$ give the four
non-metacyclic triples $(G, x, y)$ in ${\cal I}(n)$, and for $n=4$ the single non-metacyclic triple
corresponds to the group $G(4\,;0, 0)$. Let ${\cal N}(n\,;k,l)$ denote the map in ${\cal R}(n)$
corresponding to $G(n\,;k,l)$; if $k=l$ it has type $\{4,n\}$ and genus $(n-2)^2/4$, and if $k\neq l$ it has type $\{8,n\}$ and genus $1+n(3n-8)/8$. For instance, Fig.~1(b) shows ${\cal N}(4\,;0,0)$, isomorphic to  the torus map $\{4,4\}_{2,2}$ of Coxeter and Moser~\cite[Ch.~8]{CM}. The maps ${\cal N}(n\,;0,0)$ are special cases for $q=n/2=2^{e-1}$ of the maps $\{4, 2q\}_4$, the duals of the maps $\{2q, 4\}_4$ in ~\cite[Table~8]{CM}; they are also isomorphic to
the maps ${\cal O}_n$ in Example~3 of~\cite{DJKNS1}, obtained by applying Wilson's `opposite' operation~\cite{Wil} to the torus maps $\{4, 4\}_n=\{4,4\}_{q,q}$ described in~\cite[Ch.~8]{CM}. All four maps ${\cal N}(n\,;k,l)$ are $4$-sheeted regular coverings of ${\cal N}(n/2\,;0,0)$, with the covering obtained by factoring out the central subgroup $\langle x^{n/2}, y^{n/2}\rangle\cong C_2\times C_2$ of $G(n\,;k,l)$. It will be useful in \S 7 for us to define $f_G=1$ for each of these non-metacyclic $2$-groups
$G=G(n\,;k,l)$.

\section{Basic constructions}

Here we describe two basic constructions of regular embeddings of complete bipartite graphs. Our
eventual aim is to show that all such embeddings can be found by applying these constructions to the
prime power embeddings described in \S3.

\medskip

\noindent{\bf Construction~4.1\;} Let ${\cal M}_i\in{\cal R}(s_i)$ for $i=1,\ldots, k$, where the
integers $s_1,\ldots, s_k$ are mutually coprime. Each ${\cal M}_i$ corresponds to a triple $(G_i,
x_i, y_i)\in{\cal I}(s_i)$, and it is easy to check that if $G=G_1\times \cdots\times G_k$,
$x=(x_1,\ldots, x_k)$ and $y=(y_1,\ldots, y_k)$ then $(G,x,y)\in{\cal I}(s)$ where $s=s_1\ldots s_k$.
We call $(G, x, y)$ the {\sl cartesian product\/} of the triples $(G_i, x_i, y_i)$, and we call the
corresponding regular embedding of $K_{s,s}$ the {\sl cartesian product\/}
${\cal M}_1\times\cdots\times{\cal M}_k$ of the embeddings ${\cal M}_i$. For example, if $s$ has prime
power factorisation $s=p_1^{e_1}\ldots p_k^{e_k}$ then
$${\cal S}(s)\cong{\cal S}(p_1^{e_1})\times\cdots\times{\cal S}(p_k^{e_k}).$$

\smallskip

The above construction was introduced in~\cite{JNS1}; the following construction is a generalisation
of some simple examples also considered there.

\medskip

\noindent{\bf Construction~4.2\;} Let ${\cal S}\in{\cal R}(s)$ and let ${\cal T}={\cal S}(t)$,
respectively corresponding to a triple $(S, x_S, y_S)\in{\cal I}(s)$ and the standard triple $(T, x_T,
y_T)\in{\cal I}(t)$. Since $T\cong C_t\times C_t$ we can use the basis $x_Y, y_T$ of the free ${\bf
Z}_t$-module $T$ to identify ${\rm Aut}\,T$ with the general linear group $GL_2({\bf Z}_t)$. Suppose that
$s$ and $t$ are coprime, and that there is a homomorphism $\psi:S\to GL_2({\bf Z}_t)$ such that
\[x_S\mapsto
\left(\begin{array}{cc}
1 & 0 \\
0 & \lambda
\end{array}\right)
\quad{\rm and}
\quad y_S\mapsto
\left(\begin{array}{cc}
\lambda & 0\\
0 & 1
\end{array}\right)\]
for some $\lambda$ in the group ${\bf Z}^*_t$ of units in ${\bf Z}_t$; we will call this a {\sl
diagonal action\/} of $S$ on $T$, with {\sl eigenvalue\/} $\lambda$. Now let $G$ be the semidirect
product $T\negthinspace:_{\psi}\negthinspace S$ of $T$ by $S$, with $\psi$ giving the action of the
complement $S$ by conjugation on the normal subgroup $T$, so that
$$x_T^{x_S}=x_T,\quad y_T^{x_S}=y_T^{\lambda},\quad x_T^{y_S}=x_T^{\lambda}\quad{\rm and}\quad
y_T^{y_S}=y_T.$$ If we define $x=x_Sx_T$ and $y=y_Sy_T$ then these relations imply that $(G, x, y)$ is
an $n$-isobicyclic triple where $n=st$. We will call this triple a {\sl semidirect product\/} of $(T,
x_T, y_T)$ by $(S, x_S, y_S)$, and we will call the corresponding map ${\cal M}\in{\cal R}(n)$ a {\sl
semidirect product\/} ${\cal T}\negthinspace:\negthinspace{\cal S}$ of $\cal T$ by $\cal S$, writing
${\cal M}={\cal T}\negthinspace:_{\lambda}{\cal S}$ if we need to specify the value of $\lambda$. The
parameters $t$ and $\lambda$, together with the isomorphism class of $(S, x_S, y_S)$, uniquely
determine the isomorphism class of $(G, x, y)$ and hence also that of $\cal M$.  Any such map ${\cal
M}={\cal T}\negthinspace:\negthinspace{\cal S}$ is a regular covering of $\cal S$, with
$T$ as the covering group, so that ${\cal M}/T\cong{\cal S}$; we will therefore also refer to $\cal
M$ as an {\sl abelian covering\/} of $\cal S$, and similarly for the triple $(G,x,y)$.

\medskip

\noindent{\bf Example 4.3\;} Construction~4.2 yields the cartesian product ${\cal T}\times{\cal S}$ if
and only if $\lambda=1$.

\medskip

\noindent{\bf Example 4.4\;} Instances of non-trivial actions of $S$ on $T$ can be found by taking $\cal S$ to be the standard embedding, so that $S\cong C_s\times C_s$, and choosing $\lambda\in {\bf Z}^*_t$ so that $\lambda^s=1$; for example, see~\cite{JNS1} for the case where $t$ is prime, with $s$ dividing $t-1$.

\medskip

The main results of this paper are that each $n$-isobicyclic triple, and hence each regular embedding
of $K_{n,n}$, is an abelian covering of a cartesian product of prime power triples (or embeddings),
as described in \S 3. More specifically, we will prove the following:

% Theorem 4.1. 

\begin{thm} Each $n$-isobicyclic triple $(G,x,y)$ is isomorphic to a semidirect product of
a $t$-isobicyclic triple $(T, x_T, y_T)$ by an $s$-isobicyclic triple $(S, x_S, y_S)$, where
\item{\rm(i)} $n=st$ with $s$ and $t$ coprime,
\item{\rm(ii)} $(S, x_S, y_S)$ is a cartesian product of $p_i^{e_i}$-isobicyclic triples for
$i=1,\ldots, k$, where $s=p_1^{e_1}\dots p_k^{e_k}$ for distinct primes $p_i$,
\item{\rm(iii)} $(T, x_T, y_T)$ is the standard $t$-isobicyclic triple.
\end{thm}

We will prove this result in \S 8. As an immediate consequence of Theorem~4.1, using the correspondence between maps and isobicyclic triples outlined in \S 2, we have:

% Corollary 4.2. 

\begin{cor}
Each regular orientable embedding of $K_{n,n}$ has the form ${\cal M}={\cal
T}\negthinspace:\negthinspace{\cal S}$, where
\item{\rm(i)} $n=st$ with $s$ and $t$ coprime,
\item{\rm(ii)} ${\cal S}={\cal S}_1\times\cdots\times{\cal S}_k$ with each
${\cal S}_i\in{\cal R}(p_i^{e_i})$ where $s=p_1^{e_1}\dots p_k^{e_k}$ for distinct primes $p_i$,
\item{\rm(iii)} $\cal T$ is the standard embedding of $K_{t,t}$. \hfill$\square$
\end{cor}

In general, the decompositions in Theorem~4.1 and Corollary~4.2 are not unique: a Sylow subgroup which is an abelian direct factor of $G$ could appear as a direct factor of either $S$ or $T$, and there is
a corresponding result for the decomposition of $\cal M$. However, this is the only way in which
uniqueness can fail, and in each case there is a unique decomposition for which $t$ is minimal; this
is characterised by the property that $T\cap Z(G)=1$, where $Z(G)$ denotes the centre of $G$. We will
call this the {\sl canonical decomposition\/} of $(G, x, y)$ or of $\cal M$.

\section{Hall subgroups and Sylow subgroups}

If $\cal M$ is a regular embedding of $K_{n,n}$, then the corresponding group $G={\rm Aut}^+_0{\cal
M}$ is a product $XY$ of abelian subgroups $X$ and $Y$, so it follows from a theorem of It\^o~\cite{Ito}
that $G$ is metabelian, that is, an extension of one abelian group by another. Being solvable, $G$
satisfies Hall's theorems,  or the extended Sylow theorems, which we will now explain.

If $\pi$ is any set of prime numbers, then a $\pi$-{\sl number\/} is a positive integer divisible
only by primes in $\pi$, and a $\pi$-{\sl group\/} is a group whose order is a $\pi$-number. A {\sl
Hall\/} $\pi$-{\sl subgroup\/} of a group $G$ is a $\pi$-group $H\leq G$ whose index $|G:H|$ is a
$\pi'$-number, where $\pi'$ denotes the set of primes not in $\pi$. A fundamental theorem of
P.~Hall~\cite{Hal1} (see also~\cite[VI.1.8]{Hup} or~\cite[11.18]{Ros}) states that for each set $\pi$ of primes, a finite solvable group has a single conjugacy class of Hall $\pi$-subgroups, and every $\pi$-subgroup is contained in a Hall $\pi$-subgroup. For instance, if $\pi=\{p\}$ or $\{p\}'$ then the Hall $\pi$-subgroups of $G$ are simply its Sylow $p$-subgroups or Sylow $p$-complements.

A finite nilpotent group has a unique Hall $\pi$-subgroup for each set $\pi$ of primes, namely the
direct product of its Sylow $p$-subgroups for the primes $p\in\pi$. In~\cite[Satz~III.4.8]{Hup}, Huppert
has shown that if a solvable group $G$ has the form $G=XY$ for nilpotent subgroups $X$ and $Y$, and
$X_{\pi}$ and $Y_{\pi}$ are the Hall $\pi$-subgroups of $X$ and $Y$, then $X_{\pi}Y_{\pi}$ is a Hall
$\pi$-subgroup of $G$. (In fact, one can omit the hypothesis that $G$ is solvable, since this is true
for any product of two nilpotent groups, by a theorem of Kegel and Wielandt~\cite[Satz~III.4.3]{Hup}.)
This is a generalisation of results of Wielandt~\cite{Wie}, who proved these properties for Sylow
$p$-subgroups and their complements. 

If $G$ is $n$-isobicyclic, so that $G=XY$ with $X\cong Y\cong C_n$, then these results apply to $G$. We
will call $X_{\pi}Y_{\pi}$ the {\sl canonical\/} Hall $\pi$-subgroup $G_{\pi}$ of $G$. We have
$X=X_{\pi}\times X_{\pi'}$ and $Y=Y_{\pi}\times Y_{\pi'}$, so the canonical generators $x$ and $y$ of
$G$ have unique factorisations $x=x_{\pi}x_{\pi'}$ and $y=y_{\pi}y_{\pi'}$ with $x_{\pi},\,
x_{\pi'},\, y_{\pi}$ and $y_{\pi'}$ elements of $X_{\pi},\, X_{\pi'},\, Y_{\pi}$ and $Y_{\pi'}$
respectively. In fact, these four elements generate these subgroups, so we call them their {\sl
canonical generators\/}. In the particular case where $\pi=\{p\}$ for some prime $p$ dividing $n$, we
will write simply $x_p,\, x_{p'}$ and so on.

Since $X\cap Y=1$ we have $X_{\pi}\cap Y_{\pi}=1$. Moreover, the canonical automorphism $\alpha$ of
$G$, transposing $x$ and $y$, must also transpose their powers $x_{\pi}$ and $y_{\pi}$. Thus $\alpha$
leaves $G_{\pi}$ invariant and induces an automorphism $\alpha_{\pi}$ of $G_{\pi}$ transposing its
canonical generators, so $G_{\pi}$ is isobicyclic. We summarise this result as follows:

% Proposition 5.1.

\begin{prop} If $G=XY$ is $n$-isobicyclic then for each set $\pi$ of primes dividing
$n$, the group $G_{\pi}=X_{\pi}Y_{\pi}$ is an isobicyclic Hall $\pi$-subgroup of $G$. \hfill$\square$
\end{prop}

% Corollary 5.2.

\begin{cor} If $G$ is $n$-isobicyclic then for each prime $p$ dividing $n$ the
Sylow $p$-subgroups of $G$ are isomorphic to one of the isobicyclic $p$-groups described in \S3.
\end{cor}

\noindent{\sl Proof.} This follows from Proposition~5.1, the conjugacy of Sylow $p$-subgroups for
a given $p$, and the prime-power classifications summarised in \S3. \hfill$\square$

\medskip

The next result gives us an important link between the Sylow subgroups and certain normal subgroups
of $G$:

% Proposition 5.3.

\begin{prop} If $G$ is $n$-isobicyclic, $p_1,\ldots, p_l$ are the distinct primes
dividing $n$, with $p_1>\cdots>p_l$, and $P_i$ is a Sylow $p_i$-subgroup of $G$ for each $i$, then
$N_i:= P_1\ldots P_i$ is a normal subgroup of $G$ for each $i=1,\ldots, l$. 
\end{prop}

\noindent{\sl Proof.} This is a special case of a result of Wielandt~\cite[Satz~3]{Wie}, which states that
any finite group $G$ which is a product of two cyclic groups has this property, where $p_1,\ldots,
p_l$ are the primes dividing $|G|$. In our case we have $|G|=n^2$, so these are also the primes
dividing $n$. \hfill$\square$

\medskip

This result implies that $G$ has an ascending series
\[1=N_0<N_1<\cdots<N_{l-1}<N_l=G\eqno(5.1)\]
of normal (in fact, characteristic) subgroups $N_i$ of $G$ with $N_i/N_{i-1}\cong P_i$ for each
$i=1,\ldots,l$. We can (and will) without loss of generality take each $P_i$ to be the canonical
Sylow $p_i$-subgroup of $G$, so each $N_i$ is the canonical Hall $\pi$-subgroup $G_{\pi}$ of $G$ for
$\pi=\{p_1,\ldots, p_i\}$. By Corollary~5.2, each of the quotients $N_i/N_{i-1}$ in $(5.1)$ is
isomorphic to one the isobicyclic $p$-groups described in \S3, where $p=p_i$.  More generally, if $1\leq
i\leq j\leq l$ then the natural epimorphism $G\to G/N_{i-1}$ induces an isomorphism
$G_{\pi}\cong N_j/N_{i-1}$ where $\pi=\{p_i,\ldots, p_j\}$, so Proposition~5.1 implies that
$N_j/N_{i-1}$ is isobicyclic.

In order to determine the structure of $G$, and hence of $\cal M$, we need to know the canonical
Sylow subgroups $P_i$ of $G$, together with the action by conjugation of $G/N_{i-1}$ on its normal
subgroup $N_i/N_{i-1}\cong P_i$ for each $i$. (Since the quotients $N_i/N_{i-1}$ in $(5.1)$ have
mutually coprime orders, the Schur-Zassenhaus Theorem (see~\cite[I.18.1]{Hup} or~\cite [10.30]{Ros}) ensures that all the
relevant extensions split, so we do not have to consider cohomological problems associated with
possibly non-split extensions.) In considering this action, by factoring out $N_{i-1}$ if necessary,
we can assume that $i=1$, so that $N_i/N_{i-1}=N_1$ is a normal Sylow $p$-subgroup $P=P_1$ of $G$
for the largest prime $p=p_1$ dividing $n$. In this case, since $G=P_1\ldots P_l$ it is sufficient
to determine how each canonical Sylow $q$-subgroup $Q=P_j$ acts by conjugation on $P$, where $q=p_j$
for some $j>1$. We will therefore assume for the next two sections that $P$ is a normal Sylow
$p$-subgroup of $G$ for the largest prime $p=p_1$ dividing $n$, and that $Q$ is the canonical Sylow
$q$-subgroup of $G$ for a prime $q=p_j$ where $j>1$. Even in the general case, when $i>1$ and $P$ is
not necessarily normalised by $Q$, we will refer to the action of $Q$ on $P$, meaning the action
induced by the natural isomorphism of $P$ with the corresponding quotient $N_i/N_{i-1}\cong P$.

An important aid to understanding this action is given by Hall's corollary to the Burnside Basis
Theorem. Let $P$ be any finite $p$-group, and let $\Phi(P)$ denote its Frattini subgroup, the
intersection of its maximal subgroups, or equivalently the smallest normal subgroup of $P$ with an
elementary abelian quotient $\tilde P=P/\Phi(P)$. Since $\Phi(P)$ is a characteristic subgroup of $P$,
every automorphism of $P$ induces an automorphism of $\tilde P$, giving a homomorphism $\phi:{\rm
Aut}\,P\to{\rm Aut}\,\tilde P$. Hall~\cite{Hal2} (see also~\cite[III.3.18]{Hup} or~\cite[11.13(i)]{Ros}) showed that $\ker\phi$ is a $p$-group, or equivalently:

% Proposition 5.4.

\begin{prop} If $P$ is a finite $p$-group, then any $p'$-group $A$ of automorphisms of
$P$ is represented faithfully on $\tilde P:=P/\Phi(P)$.\hfill$\square$
\end{prop}

Indeed, one can regard $\tilde P$ as a $d$-dimensional vector space over the field ${\bf Z}_p$, where
$d$ is the rank (minimum number of generators) of $P$, so that $A$ acts linearly on $\tilde P$ as a
subgroup of $GL_d(p)$. In our situation, where $P$ is a normal Sylow $p$-subgroup of an isobicyclic
group $G$, we have $d=2$; any Sylow $p$-complement $K$ in $G$, acting by conjugation, induces a
$p'$-group $\overline K=K/C_K(P)$ of automorphisms of $P$, where $C_K(P)$ denotes the centraliser of
$P$ in $K$. We therefore have:

% Proposition 5.5.

\begin{prop} The group $\overline K$ is represented faithfully on $\tilde P\cong
C_p\times C_p$ as a subgroup of $GL_2(p)$. \hfill$\square$
\end{prop}

In particular, if $Q$ is a Sylow $q$-subgroup of $G$ for any prime $q\neq p$, then $\overline
Q:=Q/C_Q(P)$ is represented faithfully on $\tilde P$ as a subgroup of $GL_2(p)$.

\section{Nonabelian Sylow subgroups}

In Construction~4.2, we chose $\cal T$ to be a standard embedding. The next result explains why this is
necessary if $S$ is to act nontrivially on $T$.

% Proposition 6.1.

\begin{prop} If $G$ is isobicyclic, and $P$ is a nonabelian normal Sylow $p$-subgroup
of $G$, then $P$ is a direct factor of $G$.
\end{prop}

Note that by Construction~4.2, an abelian normal Sylow $p$-subgroup need not be a direct factor.
Returning to the series $(5.1)$ for an isobicyclic group $G$, we can immediately deduce:

% Corollary 6.2.

\begin{prop} If $P$ and $Q$ are Sylow $p$- and $q$-subgroups of an isobicyclic group $G$
for primes $p\neq q$, and $P$ is nonabelian, then $Q$ acts trivially on $P$. \hfill $\square$
\end{prop}

In order to prove Proposition~6.1, we first need the following technical result:

% Lemma 6.3.

\begin{lemma} Let $P$ be a nonabelian isobicyclic $p$-group, with canonical generators $x$ and
$y$, and canonical automorphism $\alpha$. If $A=\langle\,\beta,\,\beta^{\alpha}\,\rangle\leq{\rm
Aut}\,P$, where $\beta$ is an automorphism fixing $x$, and $A$ has order coprime to $p$, then $A=1$.
\end{lemma}

\noindent{\sl Proof.} Having order coprime to $p$, $A$ is represented faithfully on
$\tilde P:=P/\Phi(P)\cong C_p\times C_p$ by Proposition~5.4, so it is sufficient to show that $A$
acts trivially on $\tilde P$. If we use the basis of $\tilde P$ formed by the images $\tilde x$ and
$\tilde y$ of the canonical generators $x$ and $y$ of $P$, then $A$ is identified with the subgroup of
${\rm Aut}\,\tilde P=GL_2(p)$ generated by two matrices of the form
\[
M_{\beta}=\left(\begin{array}{cc}1&0\\ \kappa&\lambda \end{array}\right)\quad
{\rm and}\quad
M_{\gamma}=\left(\begin{array}{cc}\lambda&\kappa\\ 0&1\end{array}\right)
\]
with $\kappa, \lambda\in{\bf Z}_p$ and $\lambda\neq 0$. Here $M_{\beta}$ represents the action of
$\beta$, which fixes $\tilde x$, and $M_{\gamma}$ represents that of $\gamma:=\beta^{\alpha}$, which
fixes $\tilde y$. If $\lambda=1$ then such matrices have order $p$ (against our hypotheses) or $1$ (as
required), so we may assume that $\lambda\neq 1$, and hence $p>2$.

In the proof of Proposition~16 of~\cite{JNS2} it is shown that if $P$ is a non-abelian isobicyclic
$p$-group, where $p>2$, then every automorphism of $P$ acts on $\tilde P$ as a matrix in $GL_2(p)$ of
the form
\[
\left(
\begin{array}{cc}
1&*\\
0&*
\end{array}
\right)
\] with respect to the basis $\tilde g, \tilde h$ for $\tilde P$, where $g$
and $h$ are as in $(3.1)$. The group of all such matrices, isomorphic to the affine general linear group
$AGL_1(p)$, acts faithfully on the projective line $PG_1(p)$ formed by the $1$-dimensional subspaces of
$\tilde P$, fixing one element of $PG_1(p)$ and acting sharply $2$-transitively on the remaining $p$
elements. Any subgroup of $AGL_1(p)$ of order coprime to $p$ is cyclic of order dividing $p-1$, fixing
two elements and acting semi-regularly on the other $p-1$. In particular, this must apply to $A$;
however, $M_{\beta}$ and $M_{\gamma}$, between them, have at least three fixed-points on $PG_1(p)$, so $A=1$ as required.

\hfill $\square$

\medskip

\noindent{\sl Proof of Proposition 6.1.} Being normal, $P$ is the unique Sylow $p$-subgroup of $G$, so
$P$ is the canonical Sylow $p$-subgroup $G_p=X_pY_p$. Let $K$ be the canonical Sylow $p$-complement
$G_{p'}=X_{p'}Y_{p'}$ in $G$, so $G$ is a semidirect product of $P$ by $K$. Acting by conjugation on
$P$, $K$ induces a group $\overline K\cong K/C_K(P)$ of automorphisms of $P$. This is generated by
automorphisms $\beta$ and $\gamma=\beta^{\alpha}$ of $P$, induced by the canonical generators $x_{p'}$
and $y_{p'}$ of $K$. Since $X_p$ and $X_{p'}$ are contained in the abelian group $X$, $\beta$ acts
trivially on $X_p$, and similarly $\gamma$ acts trivially on $Y_p$. Since $P$ is a nonabelian isobicyclic
$p$-group, and since $\overline K$ has order coprime to $p$, Lemma~6.2 (with $A=\overline K$) implies
that $\overline K=1$, and hence $G=P\times K$ as required. \hfill $\square$

\section{Abelian Sylow subgroups}

We now assume that the normal Sylow $p$-subgroup $P$ of $G$ is abelian. Since $P$ is isobicyclic, it
follows that $P\cong C_{p^d}\times C_{p^d}$ for some $d\geq 1$. In order to determine how $G$ can act by conjugation on $P$ we need the following lemma, an analogue of Lemma~6.3:

%Lemma 7.1.

\begin{lemma} Let $P=C_{p^d}\times C_{p^d}$, and let $A\leq{\rm Aut}\,P=GL_2({\bf Z}_{p^d})$
be a $p'$-group generated by matrices
\[M_x=\left(\begin{array}{cc}1&0\\ \kappa&\lambda\end{array}\right)\quad{\it and}\quad
M_y=\left(\begin{array}{cc}\lambda&\kappa\\ 0&1\end{array}\right).\]
Then the following are equivalent:
\item{\rm(a)} $A$ is abelian,
\item{\rm(b)} $A$ consists of diagonal matrices,
\item{\rm(c)} $\kappa=0$,
\item{\rm(d)} $\kappa\equiv 0$ mod~$(p)$.
\end{lemma}

\noindent{\sl Proof.} It is clear that (b) and (c) are equivalent, that (b) implies (a),
and that (c) implies (d). To show that (a) implies (c), we have
\[M_xM_y=\left(\begin{array}{cc}\lambda&\kappa\\ \kappa\lambda&\kappa^2+\lambda\end{array}\right)\quad{\rm and}\quad
M_yM_x=\left(\begin{array}{cc}\kappa^2+\lambda&\kappa\lambda\\ \kappa&\lambda\end{array}\right),\]
so $M_x$ and $M_y$ commute if and only if $\kappa^2=0$ and $\kappa(\lambda-1)=0$. If $\kappa\neq 0$
this implies that $\lambda\equiv 1$ mod~$(p)$. Since $A$ is a $p'$-group, it is represented faithfully
on $\tilde P:=P/\Phi(P)\cong C_p\times C_p$ as a subgroup of $GL_2(p)$; however, since $\lambda\equiv
1$ mod~$(p)$ the images of $M_x$ and $M_y$ in $GL_2(p)$ have order dividing $p$, so $A=1$,
contradicting $\kappa\neq 0$. To show that (d) implies (a), if $\kappa\equiv 0$ mod~$(p)$ then the
images of $M_x$ and $M_y$ in $GL_2(p)$ are diagonal matrices, so they commute; since $A$ is
represented faithfully in $GL_2(p)$ it follows that $M_x$ and $M_y$ commute, so $A$ is abelian.
\hfill$\square$

\medskip

If $A$ satisfies the equivalent conditions in Lemma~7.1, so that
\[M_x=\left(\begin{array}{cc}1&0\\ 0&\lambda\end{array}\right)\quad{\rm and}
\quad M_y=\left(\begin{array}{cc}\lambda&0\\ 0&1\end{array}\right)\]
for some $\lambda\in{\bf Z}_{p^d}$, we will say that $A$ {\sl acts diagonally\/} on $P$ with respect to
the canonical basis of $P$, with {\sl eigenvalue\/} $\lambda$. By Proposition~5.4 we must have
$\lambda^{p-1}=1$.

% Proposition 7.2.

\begin{prop} If $P$ is an abelian normal Sylow $p$-subgroup of an isobicyclic group
$G$, then $G$ acts diagonally on $P$ with respect to its canonical basis.
\end{prop}

\noindent{\sl Proof.} Since $P$ is abelian, the action of $G$ by conjugation on $P$ is equivalent to
that of $G/P$, and this in turn is equivalent to that of the Sylow $p'$-complement $G_{p'}$ of $G$.
Since $G_{p'}$ is the product of the canonical Sylow $q$-subgroups $Q=G_q$ of $G$ for the primes
$q\neq p$ dividing $n$, it is sufficient to show that each such group $Q$ acts diagonally on $P$.

Since $|Q|$ is coprime to $p$, the subgroup $\overline Q=Q/C_Q(P)$ of ${\rm Aut}\,P$ induced by $Q$ on
$P$ is isomorphic to the subgroup $\tilde Q$ of $GL_2(p)$ induced by $Q$ on $\tilde P:=P/\Phi(P)$. Now
$\tilde Q$ is contained in a Sylow $q$-subgroup $S$ of $GL_2(p)$. If $q$ is odd then $S$ is abelian,
conjugate to a subgroup of the diagonal group $GL_1(p)\times GL_1(p)$ if $q$ divides $p-1$, or to a
subgroup of a Singer group $GF(p^2)^*\cong C_{p^2-1}$ if not. In either case, $\tilde Q$ is abelian
and hence so is $\overline Q$, so Lemma~7.1 implies that $Q$ acts diagonally on $P$.

We may therefore assume that $q=2$. The above argument does not apply in this case, since a Sylow
$2$-subgroup $S$ of $GL_2(p)$ is nonabelian for each prime $p>2$. If $p\equiv 1$ mod~$(4)$ then $S$
is a wreath product $C_{2^r}\wr C_2$ where $2^r\parallel p-1$: one can take the base group
$C_{2^r}\times C_{2^r}$ to consist of diagonal matrices, with the complement $C_2$ generated by
\(
\left(
\begin{array}{cc}
0&1\\
1&0
\end{array}
\right).
\)
If $p\equiv -1$ mod~$(4)$  then $S$ is a semidihedral group
\[SD_{2^{r+1}}=\langle c,d\mid c^{2^{r+1}}=d^2=1, c^d=c^{-1+2^r}\rangle\]
where $2^r\parallel p+1$: the normal subgroup $\langle c\rangle\cong C_{2^{r+1}}$ is contained in a
Singer group $GF(p^2)^*$, and the complement $\langle d\rangle\cong C_2$ is induced by the Galois
group of $GF(p^2)$. In either case, $S$ has an abelian subgroup of index $2$, so all squares in $Q$
commute.

The canonical generators of $Q$, acting by conjugation on $P$, are represented by matrices of the form
$M_x, M_y\in GL_2({\bf Z}_{p^d})$ in Lemma~7.1. Now
\[M_x^2=\left(\begin{array}{cc}1&0\\ \kappa(\lambda+1)&\lambda^2\end{array}\right)\quad{\rm and}\quad
M_y^2=\left(\begin{array}{cc}\lambda^2&\kappa(\lambda+1)\\ 0&1\end{array}\right),\]
so
\[M_x^2M_y^2=\left(\begin{array}{cc}\lambda^2&\kappa(\lambda+1)\\
\kappa\lambda^2(\lambda+1)&\kappa^2(\lambda+1)^2+\lambda^2\end{array}\right)\]
and
\[M_y^2M_x^2=\left(\begin{array}{cc}\kappa^2(\lambda+1)^2+\lambda^2&\kappa \lambda^2(\lambda+1)\\
\kappa(\lambda+1)&\lambda^2\end{array}\right).\]
This shows that $M_x^2$ and $M_y^2$ commute if and only if $\kappa^2(\lambda+1)^2=0$ and
$\kappa(\lambda+1)(\lambda^2-1)=0$ in ${\bf Z}_{p^d}$. Since they generate a $p'$-group, the matrices
$M_x^2$ and $M_y^2$ commute if and only if their images in $GL_2(p)$ commute, that is, 
$\kappa^2(\lambda+1)^2\equiv \kappa(\lambda+1)(\lambda^2-1)\equiv 0$ mod~$(p)$, or equivalently
$\kappa\equiv 0$ or $\lambda\equiv -1$ mod~$(p)$. If $\kappa\equiv 0$ mod~$(p)$ then Lemma~7.1 implies
that $Q$ acts diagonally, so assume that $\lambda\equiv -1$ mod~$(p)$. Then $M_x^2$ and
$M_y^2$ induce the identity on $\tilde P$, so they also induce the identity on $P$. Thus the images
$\overline{X_q}=\langle M_x\rangle$ and $\overline{Y_q}=\langle M_y\rangle$ of $X_q$ and $Y_q$ in
${\rm Aut}\,P$ each have order at most $2$, so $\overline Q=\overline{X_q}\overline{Y_q}$ has order
at most $4$. It follows that $\overline Q$ is abelian, so Lemma~7.1 implies that $Q$ acts diagonally
on $P$.
\hfill$\square$

\medskip

We can now deduce the analogue of Corollary~6.2 for abelian Sylow $p$-subgroups:

% Corollary 7.3.

\begin{cor} If $P$ and $Q$ are Sylow $p$- and $q$-subgroups of an isobicyclic group $G$
for primes $p\neq q$, and $P$ is abelian, then $Q$ acts diagonally on $P$. If $q$ does not divide
$p-1$ then this action is trivial.
\end{cor}

\noindent{\sl Proof.} The first statement follows immediately from Proposition~7.2.
Proposition~5.5 then implies that $\overline Q$ is represented faithfully on $\tilde P$ as a subgroup
of $GL_1(p)\times GL_1(p)\cong C_{p-1}\times C_{p-1}$, so if $q$ does not divide $p-1$ then $\overline
Q=1$, that is, $Q$ acts trivially on $P$.
\hfill$\square$

\medskip

It remains for us to determine how $Q$ can act diagonally on $P$ when $q$ divides $p-1$.

% Lemma 7.4.

\begin{lemma} Let $p$ and $q$ be primes, with $q$ dividing $p-1$. An isobicyclic $q$-group $Q$ can
act diagonally on
$P\cong C_{p^d}\times C_{p^d}$ with eigenvalue $\lambda\in{\bf Z}_{p^d}$ if and only if
$\lambda^{q^m}=1$, where
$m=\min\{f_Q,r\}$ and $q^r\parallel p-1$; there are $q^m$ such elements $\lambda\in {\bf Z}_{p^d}$.
\end{lemma}

\noindent{\sl Proof.} Since diagonal matrices commute, there is a diagonal action of $Q$ on $P$ with
eigenvalue $\lambda$ if and only if the mapping
\[x\mapsto M_x=\left(\begin{array}{cc}1&0\\ 0&\lambda\end{array}\right)\quad
{\rm and}\quad y\mapsto M_y=\left(\begin{array}{cc}\lambda&0\\ 0&1\end{array}\right)\]
of the canonical generators $x$ and $y$ of $Q$ induces a homomorphism $Q^{\rm ab}:=Q/Q'\to GL_2({\bf Z}_{p^d})$.

First suppose that $Q$ is metacyclic. The results in~\cite{JNS2} and~\cite{DJKNS1} (see \S3) show that
$$Q=\langle g,h\mid g^{q^e}=h^{q^e}=1,\, h^g=h^{1+q^f}\rangle$$
where $|Q|=q^{2e}$, with $f=f_Q$ taking a value $1,\ldots, e$ if $q>2$, or $2,\ldots, e$ if $q=2$, so
\[Q^{\rm ab}=\langle g, h\mid g^{q^e}=h^{q^f}=[g,h]=1\rangle.\]
The canonical generators of $Q$ are $x=g^u$ and $y=g^uh$ where $u=1,\ldots, q^{e-f}$ is coprime to $q$, so we have an equivalent presentation
\[Q^{\rm ab}=\langle x,y\mid x^{q^e}=[x,y]=1,\, x^{q^f}=y^{q^f}\rangle.\]
Since $f\leq e$ we see that $M_x$ and $M_y$ satisfy these relations, giving a diagonal action of $Q$
with eigenvalue $\lambda$, if and only if $\lambda^{q^f}=1$. Since $q$ divides $p-1$ we have $p>2$, so the multiplicative group ${\bf Z}_{p^d}^*$ is cyclic, of order $\phi(p^d)=(p-1)p^{d-1}$. The number of solutions of $\lambda^{q^f}=1$ in this group is therefore $\gcd(q^f,\phi(p^d))=\gcd(q^f,p-1)=\gcd(q^f,q^r)=q^m$ where $q^r\parallel p-1$ and $m=\min\{f,r\}$.

Now suppose that $Q$ is not metacyclic, so $q=2$. As shown in~\cite{DJKNS2} and explained in \S3, we have
\begin{eqnarray*}
Q=\langle x,y\; \mid & x^{2^e}=y^{2^e}=1,\,c:=[y,x]=x^{2+k2^{e-1}}y^{-2-k2^{e-1}},\\
& c^x=c^{-1+l2^{e-2}}x^4,\, c^y=c^{-1-l2^{e-2}}y^{-4}\rangle,
\end{eqnarray*}
where $e\geq 2$ and $k,l\in\{0,1\}$, with $k=l=0$ if $e=2$. Thus
\[Q^{\rm ab}=\langle x,y\mid x^4=y^4=[x,y]=1,\, x^{2+k2^{e-1}}=y^{2+k2^{e-1}}\rangle,\]
so we obtain a diagonal action of $Q$ if and only if $\lambda^2=1$. Since $p$ is odd it follows that $Q$
has two diagonal actions on $P$, given by $\lambda=\pm 1$. Since we defined $f_Q=1$ for the
non-metacyclic groups $Q$, so that $m=1$ and hence $q^m=2$, this is consistent with the statement of the
Lemma.
\hfill$\square$

% Proposition 7.5.

\begin{prop} Let $P, Q$ and $R$ be Sylow $p$-, $q$- and $r$-subgroups of an
$n$-isobicyclic group $G$, for primes $p, q$ and $r$ with $p>q>r$. If $R$ acts nontrivially on $Q$,
then $Q$ acts trivially on $P$.
\end{prop}

\noindent{\sl Proof.} Without loss of generality we can take $P, Q$ and $R$ to be the canonical
Sylow subgroups $G_p, G_q$ and $G_r$ of $G$. Since $R$ acts nontrivially on $Q$, it follows from
Proposition~6.1 (applied to the action by conjugation on $Q$, rather than $P$) that $Q$ must be
abelian, so $Q\cong C_{q^e}\times C_{q^e}$ where $q^e\parallel n$. It then follows from
Proposition~7.2 (again applied to $Q$ rather than $P$) that $R$ acts diagonally on $Q$, with an
eigenvalue $\lambda\neq 1$ in ${\bf Z}_{q^e}$. If $x_q$ and $y_q$ are the canonical generators
of $Q$, and $x_r$ and $y_r$ are those of $R$, then $x_q^{y_r}=x_q^{\lambda}$, so
$[x_q,y_r]=x_q^{\lambda-1}$, and similarly $[y_q, x_r]=y_q^{\lambda-1}$. Now $\lambda\not\equiv 1$
mod~$(q)$, since the $q'$-group $\langle \lambda\rangle\leq {\bf Z}^*_{q^e}$, acting by multiplication on
the additive $q$-group ${\bf Z}_{q^e}$, is represented faithfully on ${\bf Z}_{q^e}/\Phi({\bf
Z}_{q^e})={\bf Z}_{q^e}/q{\bf Z}_{q^e}\cong {\bf Z}_q$. It follows that $x_q^{\lambda-1}$ and
$y_q^{\lambda-1}$ are generators of the groups $\langle x_q\rangle=X_q$ and $\langle y_q\rangle=Y_q$,
so the subgroup $Q=X_qY_q$, generated by two commutators, is contained in the derived group $G'$ of
$G$. If $P$ is abelian or nonabelian, then $G$ acts diagonally or trivially on $P$, by Propositions~7.2
and 6.1 respectively. In either case, $G$ induces an abelian group of automorphisms of $P$, and since
$Q\leq G'$ it follows that $Q$ acts trivially on $P$. \hfill$\square$

\section{Isobicyclic groups and directed graphs}

Using Corollary~7.3 and Lemma~7.5, we can describe the structure of an isobicyclic triple $(G,
x, y)\in{\cal I}(n)$ by means of a labelled directed graph. We will use this idea to give a proof of
Theorem~4.1.

Let $\Pi$ denote the set of all prime numbers, and let $\to$ denote the binary relation on $\Pi$
defined by $q\to p$ if and only if $q$ divides $p-1$; we can regard $\Pi$ as a directed graph, with
an arc from $q$ to $p$ whenever $q\to p$. For each integer $n\geq 2$, let $\Pi_n$ denote the induced
subgraph of $\Pi$ whose vertices are the prime factors $p_1,\ldots, p_l$ of $n$, formed by
restricting the relation $\to$ to these primes. We say that $\rup$ is a {\sl short subrelation\/} of $\to$
on $\Pi_n$ if $\rup$ is a subrelation of $\to$ on $\Pi_n$ (i.e.~$q\rup p$ implies that $q\to p$) and there
is no triple $p, q, r$ in $\Pi_n$ with $r\rup q\rup p$; equivalently, $\rup$ defines a directed subgraph
$\Gamma$ of $\Pi_n$ which spans $\Pi_n$ (i.e.~contains all the vertices of $\Pi_n$), and contains no
directed paths of length greater than $1$. In these circumstances we call $\Gamma$ a {\sl short spanning
subgraph\/} of $\Pi_n$, written $\Gamma\preceq\Pi_n$.

The motivation for these definitions is that each short subrelation $\rup$ (equivalently each short
spanning subgraph $\Gamma$) represents a choice of which Sylow subgroups of an $n$-isobicyclic group $G$ act nontrivially by conjugation on each other, with $q\rup p$ meaning that a Sylow $q$-subgroup $Q$ acts nontrivially on a Sylow $p$-subgroup $P$. By Corollary~7.3, if $q\rup p$ then $q$ divides $p-1$, so
$\rup$ must be a subrelation of $\to$. Proposition~7.5 implies that there can be no directed paths
$r\rup q\rup p$ in $\Gamma$, so this subrelation must be short.

In order to specify an $n$-isobicyclic triple $(G, x, y)$ completely, we need to describe the triples
$(P, x_p, y_p)$ corresponding to the Sylow subgroups $P=G_p$ of $G$, together with the actions
$Q\to{\rm Aut}\,P$ of these subgroups by conjugation on each other. We do this by attaching labels to
the vertices and arcs of $\Gamma$, describing these triples and actions. First we define
$$T(\Gamma)=\{p\in\Pi_n\mid q\rup p\;\;\hbox{for some}\;\;q\in\Pi_n\},$$
the set of {\sl terminal\/} vertices of $\Gamma$, and
$$N(\Gamma)=\{p\in\Pi_n\mid q\rup p\;\;\hbox{for no}\;\;q\in\Pi_n\},$$
the set of {\sl non-terminal\/} vertices. An {\sl isobicyclic labelling\/} $\Lambda$ of $\Gamma$ has
two ingredients. Firstly, each vertex $p$ of $\Gamma$ is labelled with a triple
$\Lambda(p)=(P, x_p, y_p)\in{\cal I}(p^d)$ where $p^d\parallel n$, with the restriction that if $p\in
T(\Gamma)$ then $\Lambda(p)$ is the standard triple in ${\cal I}(p^d)$. Secondly, each arc $q\rup p$ in
$\Gamma$ is labelled with an element $\lambda=\Lambda(q\rup p)\in{\bf Z}_{p^d}$ satisfying
$\lambda^{q^m}=1\neq\lambda$, where $q$ is labelled with $\Lambda(q)=(Q, x_q, y_q)$, $q^r\parallel p-1$, and $m=\min\{f_Q,r\}$ (see \S3 for the definition of $f_Q$).

Each $n$-isobicyclic triple $(G, x, y)$ determines an isobicyclic labelling $\Lambda=\Lambda(G,x,y)$
of a short spanning subgraph $\Gamma=\Gamma(G, x, y)\preceq\Pi_n$: the triple $\Lambda(p)=(P, x_p, y_p)$ labelling a vertex $p$ of $\Pi_n$ represents the canonical Sylow $p$-subgroup of $G$ and its canonical generators, with Corollary~5.2 implying that this triple is isobicyclic; an arc $q\rup p$ in $\Gamma$ corresponds to a nontrivial action $Q\to{\rm Aut}\,P$ between Sylow subgroups of $G$, with
Corollary~6.2 implying that $P$ must then be abelian, so that $\Lambda(p)$ is the standard triple;
these arcs define a short subgraph $\Gamma\preceq\Pi_n$ by Proposition~7.5 and Corollary~7.3; if
$\Lambda(q)=(Q, x_q, y_q)$ then the arc $q\rup p$ is labelled with the common non-identity eigenvalue
$\lambda$ of $x_q$ and $y_q$ in the diagonal action of $Q$ on $P$ (see Corollary~7.3), and this
satisfies $\lambda^{q^m}=1\neq\lambda$ by Lemma~7.4.

Let ${\cal L}(n)$ denote the set of all pairs $(\Gamma,\Lambda)$ where $\Lambda$ is an isobicyclic
labelling of a short spanning subgraph $\Gamma$ of $\Pi_n$. We will often refer to these simply as
{\sl isobicyclic labellings\/}. Let ${\cal J}(n)$ denote the subset of ${\cal I}(n)$ consisting of the
triples $(G, x, y)$ satisfying the conclusions of Theorem~4.1. Proving Theorem~4.1 then amounts to showing that ${\cal J}(n)={\cal I}(n)$.

% Proposition 8.1.

\begin{prop} Each isobicyclic labelling $(\Gamma,\Lambda)\in{\cal L}(n)$ is induced
by a unique triple $(G, x, y)\in{\cal I}(n)$; this triple is in ${\cal J}(n)$.
\end{prop}

\noindent{\sl Proof.} We first prove the existence of a suitable triple $(G, x, y)\in{\cal J}(n)$ by
induction on the number $l=|\Pi_n|$ of primes dividing $n$. The case $l=1$ is trivial, with $(G,
x, y)=\Lambda(p)\in{\cal J}(n)$ for the unique prime $p$ dividing $n$, so assume that $l>1$ and that
existence has been proved for all integers divisible by $l-1$ primes.

Let $p$ be the largest prime dividing $n$, so there are no arcs $p\rup q$ in $\Gamma$, and let
$\Lambda(p)=(P, x_P, y_P)$. Let $\Gamma^*$ be the directed graph formed from $\Gamma$ by deleting
the vertex $p$ and any incident arcs $q\rup p$, and let $\Lambda^*$ be the restriction of the
labelling $\Lambda$ to $\Gamma^*$. Then $(\Gamma^*,\Lambda^*)\in{\cal L}(n^*)$ where $n^*=n/p^d$ with $p^d\parallel n$, so by the induction hypothesis this labelling corresponds to a triple $(G^*, x^*,
y^*)=(T^*, x_{T^*}, y_{T^*}):(S^*, x_{S^*}, y_{S^*})\in{\cal J}(n^*)$. Without loss of generality, we can
assume that this is the canonical decomposition of this triple, so $S^*$ and $T^*$ are the direct
products of the Sylow $q$-subgroups $Q$ of $G^*$ corresponding to the non-terminal and terminal
vertices $q$ of $\Gamma^*$.

If $P$ is nonabelian then there are no arcs $q\rup p$ in $\Gamma$, so the triple $(G, x, y)=(P, x_P,
y_P)\times(G^*, x^*, y^*)\in{\cal I}(n)$ induces the labelling $(\Gamma,\Lambda)$. This triple has the
form $(T, x_T, y_T):(S, x_S, y_S)$ with $(T, x_T, y_T)=(T^*, x_{T^*}, y_{T^*})$ and $(S, x_S, y_S)=(P,
x_P, y_P)\times(S^*, x_{S^*}, y_{S^*})$, so it is in ${\cal J}(n)$.

If $P$ is abelian then $P\cong C_{p^d}\times C_{p^d}$. We need to show that there is a diagonal
action of $G^*$ on $P$ which is consistent with the labelling in $\Lambda$ of the arcs $q\rup p$ of
$\Gamma$. No terminal vertex $q$ of $\Gamma^*$ can be the source of an arc $q\rup p$, since
$\Gamma$ is short, so $T^*$ must be in the kernel of such an action. We therefore define the action
$G^*\to{\rm Aut}\,P$ by composing the natural epimorphism $G^*\to G^*/T^*\cong S^*$ with the unique
extension to $S^*$ of the actions $Q\to{\rm Aut}\,P$ for the Sylow $q$-subgroups $Q$ of $G^*$
corresponding to the  non-terminal vertices $q$ of $\Gamma^*$. This extension exists since the direct
factors $Q$ of $S^*$ commute with each other and are represented by diagonal (and hence commuting)
matrices on $P$; it is unique since these subgroups $Q$ generate $S^*$. (Indeed, the eigenvalue $\lambda$ for this action of $G^*$ is the product of the labels $\Lambda(q\rup p)$ on the arcs $q\rup p$ in
$\Gamma$.) We now form the extension $(G, x, y)=(P, x_P, y_P):_{\lambda}(G^*, x^*, y^*)$ corresponding to this action. By its construction, this triple is in ${\cal I}(n)$ and it induces the labelling
$(\Gamma,\Lambda)$. To see that $(G, x, y)\in{\cal J}(n)$, note that the abelian normal groups $P$
and $T^*$ commute and have coprime orders, so they generate an abelian normal subgroup $T=P\times T^*$ in $G$; this is complemented by $S^*$, which acts diagonally on $T$, so $(G, x, y)=(T, x_T, y_T):(S,
x_S, y_S)$ with $(T, x_T, y_T)=(P, x_P, y_P)\times(T^*, x_{T^*}, y_{T^*})$ and $(S, x_S, y_S)=(S^*,
x_{S^*}, y_{S^*})$.

The uniqueness of $(G, x, y)$ is similarly proved by induction on $l$: any triple in ${\cal
I}(n)$ with labelling $(\Gamma,\Lambda)$ must have the form $(P, x_P, y_P)\times(G^*, x^*, y^*)$ or $(P,
x_P, y_P):_{\lambda}(G^*, x^*, y^*)$ as above; such a product is determined up to isomorphism by the two
factors (and, in the latter case, the eigenvalue $\lambda$), and these in turn are uniquely determined by
$(\Gamma,\Lambda)$. \hfill$\square$

\medskip

\noindent{\sl Proof of Theorem 4.1.} Any $(G, x, y)\in{\cal I}(n)$ induces an isobicyclic labelling
$(\Gamma,\Lambda)\in{\cal L}(n)$. By Proposition~8.1 this is induced by a unique isobicyclic triple,
which is in ${\cal J}(n)$. Thus $(G, x, y)$ is in ${\cal J}(n)$, so it satisfies the conclusions
of Theorem~4.1. \hfill$\square$

\medskip

Corollary~4.2, which gives a similar decomposition for the maps ${\cal M}\in{\cal R}(n)$, is an
immediate consequence of Theorem~4.1.

As observed at the end of \S 4, these decompositions of triples and maps are not in general unique.
We will now discuss this lack of uniqueness in terms of isobicyclic labellings; we will do so mainly for
triples, but the analogous comments apply in the obvious way to maps. Each triple
$(G,x,y)\in{\cal I}(n)$ corresponds to a labelling $(\Gamma,\Lambda)\in{\cal L}(n)$. This triple can be
decomposed as a cartesian product of triples in ${\cal I}(n')$ for various divisors $n'$ of $n$, with the
indecomposable direct factors corresponding to the connected components of $\Gamma$, each with its
labelling induced by $\Lambda$. In particular, any isolated vertex $p$ of $\Gamma$ represents a Sylow
$p$-subgroup $P$ which is a direct factor of $G$, and if $\Lambda(p)$ is a standard triple then the
abelian group $P$ can be regarded as a subgroup of $S$ or $T$. This shows that the decomposition in
Theorem~4.1 is not generally unique. However, any terminal vertex $p$ must correspond to a Sylow
$p$-subgroup $P\leq T$, so if we use the canonical decomposition, with $t$ minimal, then $T$ is the
Hall subgroup $G_{\pi}$ for $\pi=T(\Gamma)$, the direct product of these Sylow subgroups $P$ as in the
proof of Proposition~8.1. We will assume from now on that the canonical decomposition is always used,
unless stated otherwise. Then $\cal M$, or equivalently $(G, x, y)$, uniquely determines $t$ and the
direct factors ${\cal S}_i$ of $\cal S$ in Corollary~4.2, up to a permutation of these factors. It also
uniquely determines the eigenvalue $\lambda$ for the diagonal action of $S$ on $T$; this is shown by
the following result, which also establishes some properties of $\lambda$ needed in the next section:

% Lemma 8.2.

\begin{lemma}
In the canonical decompositions for $\cal M$ and $(G, x, y)$ we have:
\item{\rm(i)} The eigenvalue $\lambda$ for $S$ on $T$ is the unique solution in ${\bf Z}_t$ of
the congruences $\lambda\equiv \lambda_p$ mod~$(p^d)$ for each $p\in\pi=T(\Gamma)$, where
$p^d\parallel n$ and $\lambda_p$ is the product in ${\bf Z}_{p^d}^*$ of the labels $\Lambda(q\rup p)$
attached to the arcs $q\rup p$ in $\Gamma$.
\item{\rm(ii)} The multiplicative order $|\lambda|$ of $\lambda$ in ${\bf Z}_t^*$ is ${\rm
lcm}\{|\lambda_p|\mid p\in\pi\}$, where the order $|\lambda_p|$ of $\lambda_p$ in ${\bf Z}_{p^d}^*$
is the product of the orders of the labels $\Lambda(q\rup p)\in{\bf Z}_{p^d}^*$ on the arcs $q\rup p$
in $\Gamma$.
\item{\rm(iii)} If $t>1$ then $\lambda-1$ is a unit in ${\bf Z}_t$.
\end{lemma}

\noindent{\sl Proof.} (i) $S$ acts diagonally on the Sylow $p$-subgroup $P$ of $T$ for each prime
$p\in\pi$. Since $x_S=x_{\pi'}=\prod_{q\in\pi'}x_q$, and similarly for $y_S$, the eigenvalue
$\lambda_p$ for this action of $S$ on $P$ is the product of the eigenvalues for its Sylow
$q$-subgroups, and this is the product of the labels $\Lambda(q\rup p)$ on the arcs $q\rup p$ in
$\Gamma$. By the Chinese Remainder Theorem, there is a unique $\lambda\in{\bf Z}_t^*$ such that
$\lambda\equiv \lambda_p$ mod~$(p^d)$ for each $p\in\pi$, and since $T$ is the direct sum of its Sylow
$p$-subgroups, this is the eigenvalue for $S$ on $T$.

\smallskip

\noindent(ii) For a given $p$, the eigenvalues $\Lambda(q\rup p)$ have mutually coprime orders (since
each is a power of $q$ by Lemma~7.4), so the multiplicative order $|\lambda_p|$ of $\lambda_p$ in
${\bf Z}_{p^d}^*$ is the product of the orders of these eigenvalues. Since $\lambda$ is identified
with $(\lambda_p)$ under the natural isomorphism ${\bf Z}_t^*\cong\prod_{p\in\pi}{\bf Z}_{p^d}^*$, it
has multiplicative order $|\lambda|={\rm lcm}\{|\lambda_p|\mid p\in\pi\}$ in ${\bf Z}_t^*$.

\smallskip

\noindent{\rm(iii)} By the minimality of $t$, each $p\in\pi$ is a terminal vertex of $\Gamma$, so
$\lambda_p\neq 1$. Then Proposition~5.4 implies that $\lambda_p-1$ is a unit in ${\bf Z}_{p^d}$
for each $p$, and since $\lambda\equiv \lambda_p$ mod~$(p^d)$ for each $p$ it follows that
$\lambda-1$ is a unit in ${\bf Z}_t$.
\hfill$\square$

\medskip

Note that the exceptional case $t=1$ in (iii) corresponds to the situation where ${\cal M}={\cal
S}={\cal S}_1\times\cdots\times{\cal S}_k$ in Corollary~4.2, with $\lambda=1$.

\section{Type and genus}

In the notation of~\cite{CM}, a regular map has {\sl type\/} $\{p, q\}$ if its faces are $p$-gons and its
vertices have valency $q$. If $\cal M$ is a regular embedding of $K_{n,n}$, corresponding to a triple
$(G, x, y)\in{\cal I}(n)$, then the faces of $\cal M$ are all $2m$-gons, so $\cal M$ has type
$\{2m,n\}$, where $m=m_G$ is the order $|xy|$ of $xy$ in $G$. It follows that there are $n^2/m$ faces, so
$\cal M$ has Euler characteristic
\[\chi=2n-n^2+\frac{n^2}{m}\]
and genus
\[g=1-\frac{\chi}{2}=1+\frac{n}{2}\Bigl(n-\frac{n}{m}-2\Bigr).\]
For instance, it is shown in~\cite{JNS2} and~\cite{DJKNS1} that if $\cal M$ is a metacyclic prime power
embedding, then $m=n$ and so $\cal M$ has type $\{2n,n\}$ and genus $(n-1)(n-2)/2$. However,
Lemma~4.3(4) of~\cite{DJKNS2} states that for $n=2^e$ the non-metacyclic $2$-groups $G(n\,;k,l)$ described
in \S3 satisfy
\[y^jx^i=x^{-i}y^{-j}z^{k+l({i+\frac{j}{2}})}\]
for all odd $i$ and $j$, where $z$ is the central involution $x^{n/2}y^{n/2}$, so $(xy)^2=z^{k+l}$; it
follows that the non-metacyclic embeddings ${\cal N}(n\,;0,0)$ and ${\cal N}(n\,;1,1)$ have
$m=2$, so they have type $\{4,n\}$ and genus $(n-2)^2/4$, while ${\cal N}(n\,;0,1)$ and ${\cal
N}(n\,;1,0)$ have $m=4$, so they have type $\{8 ,n\}$ and genus $1+n(3n-8)/8$.

% Lemma 9.1.

\begin{lemma}
If $S$ is as in the canonical decomposition of $\cal M$ then $m_G=m_S$.
\end{lemma}

\noindent{\sl Proof.} By Corollary~4.2 we have $x=x_Sx_T$ and $y=y_Sy_T$, so
\[xy=x_Sx_T.y_Sy_T=x_Sy_S.y_S^{-1}x_Ty_Sy_T=ab\]
where
\[a=x_Sy_S\in S\quad{\rm and}\quad b=y_S^{-1}x_Ty_Sy_T\in T.\]
We can identify $T$ with ${\bf Z}_t^2$ so that $x_T$ and $y_T$ correspond to the basis elements
$(1,0)$ and $(0,1)$. Since $y_S$ acts by conjugation as the matrix
\( \left( \begin{array}{cc} \lambda & 0 \\ 0 & 1 \end{array} \right) \)
with respect to this basis, we have $y_S^{-1}x_Ty_S=(\lambda,0)$ and so $b=(\lambda,1)$. Since $x_S$ acts as
\( \left( \begin{array}{cc} 1 & 0 \\ 0 &\lambda \end{array} \right) \)
, the element $a=x_Sy_S$ acts as $\lambda I$. For any $i\geq 1$ we have
\[(ab)^i=a^i.b^{a^{i-1}}b^{a^{i-2}}\negthinspace\ldots b^ab\]
with $a^i\in S$ and each $b^{a^j}\in T$; since $S\cap T=1$ it follows that
\[(xy)^i=1 \iff a^i=1\quad{\rm and}\quad b^{a^{i-1}}b^{a^{i-2}}\negthinspace\ldots b^ab=1.\]

We have $a^i=1$ if and only if $i$ is divisible by the order $m_S$ of $x_Sy_S$ in $S$. Each $b^{a^j}$
is identified with $\lambda^j(\lambda,1)=(\lambda^{j+1},\lambda^j)$ in ${\bf Z}_t^2$, so
\[b^{a^{i-1}}b^{a^{i-2}}\negthinspace\ldots
b^ab=\sum_{j=0}^{i-1}(\lambda^{j+1},\lambda^j)=(\lambda+\lambda^2+\cdots+\lambda^i,
1+\lambda+\cdots+\lambda^{i-1}).\]
Thus $b^{a^{i-1}}b^{a^{i-2}}\negthinspace\ldots b^ab=1$ if and only if
$1+\lambda+\cdots+\lambda^{i-1}=0$ in ${\bf Z}_t$, and since $\lambda-1$ is a unit by Lemma~8.2(iii),
this is equivalent to $\lambda^i=1$, that is, to $i$ being divisible by the multiplicative order
$|\lambda|$ of $\lambda$ in ${\bf Z}_t^*$. Since $\lambda$ is an eigenvalue of $x_Sy_S$ on $T$,
$|\lambda|$ divides $m_S$, so $(xy)^i=1$ if and only if $i$ is divisible by $m_S$. Thus $xy$ has
order $m_S$. \hfill $\square$
 
\medskip

Since $S$ is the direct product of its Sylow $q$-subgroups $Q$, we have $m_S=\prod_{q\in\pi'}m_Q$
where $m_Q$ is the order of $x_qy_q$ in $Q$ and $\pi=T(\Gamma)$. For example, if $Q$ is metacyclic (as
must happen if
$q>2$) then $m_q=q^e$ where $q^e\parallel n$; it follows that $m_S=s=n_{\pi'}$, the maximal $\pi'$-number
dividing $n$, unless $S$ has a non-metacyclic Sylow $2$-subgroup, in which case $m_S=2n_{\pi^*}$ or
$4n_{\pi^*}$ where $\pi^*=\pi'\setminus\{2\}$. We therefore have:

% Proposition 9.2.

\begin{prop}
Let ${\cal M}\in{\cal R}(n)$, corresponding to a pair $(\Gamma,
\Lambda)\in{\cal L}(n)$, and let $\pi$ be the set of primes which are terminal vertices of
$\Gamma$.
\item{\rm(i)} If $n$ is odd, or if $n$ is even and the Sylow $2$-subgroup $\Lambda(2)$ is metacyclic,
then $\cal M$ has type $\{2n_{\pi'},n\}$ and genus
\[1+\frac{n}{2}\Bigl(n-n_{\pi}-2\Bigr).\]
\item{\rm(ii)} If $n$ is even and $\Lambda(2)$ is $G(n\,;0,0)$ or $G(n\,;1,1)$ then $\cal M$ has
type $\{4n_{\pi^*},n\}$ and genus
\[1+\frac{n}{2}\Bigl(n-\frac{n_{\pi^+}}{2}-2\Bigr)\]
where $\pi^+=\pi\cup\{2\}$.
\item{\rm(iii)} If $n$ is even and $\Lambda(2)$ is $G(n\,;1,0)$ or $G(n\,;0,1)$ then $\cal M$ has
type $\{8n_{\pi^*},n\}$ and genus
\[1+\frac{n}{2}\Bigl(n-\frac{n_{\pi^+}}{4}-2\Bigr)\eqno{\square}\]
\end{prop}

\section{Operations on maps}

Wilson's operations $H_j$, introduced in~\cite{Wil}, act on maps $\cal M$ by raising the rotation of edges
around each vertex to its $j$-th power, where $j$ is coprime to the valencies; they preserve the
embedded graph, the orientability and the automorphism group of $\cal M$. For instance,
$H_{-1}({\cal M})$ is the mirror image $\overline{\cal M}$ of $\cal M$. If a map ${\cal M}\in{\cal
R}(n)$ corresponds to an isobicyclic triple $(G, x, y)\in{\cal I}(n)$, and $j$ is coprime to $n$, then
$H_j({\cal M})$ is the map in ${\cal R}(n)$ corresponding to the triple $H_j(G, x, y)=(G, x^j,
y^j)\in{\cal I}(n)$. Similarly, if $\cal M$ corresponds to a labelling $(\Gamma,\Lambda)\in{\cal L}(n)$
then $H_j({\cal M})$ corresponds to the labelling $(\Gamma,\Lambda^j)\in{\cal L}(n)$ defined by
$\Lambda^j(p)=H_j(\Lambda(p))$ for each vertex $p$ of $\Gamma$ and $\Lambda^j(q\rup p)=\Lambda(q\rup
p)^j$ for each arc $q\rup p$. This gives isomorphic actions of the multiplicative group ${\bf Z}_n^*$ on
the sets ${\cal R}(n)$, ${\cal I}(n)$ and ${\cal L}(n)$. For instance, if $\Lambda(p)$ is a metacyclic
triple $(P, x_p, y_p)$, corresponding to a map ${\cal M}(p^e, f, u)\in{\cal R}(p^e)$ in the notation of
\S3, then $\Lambda^j(p)$ corresponds to the map $H_j({\cal M}(p^e, f, u))={\cal M}(p^e, f, ju)\in{\cal
R}(p^e)$; for a given pair $p^e$ and $f$ (equivalently, for a given $p^e$-isobicyclic group $P=G_f$),
these maps form a single orbit of ${\bf Z}_n^*$, as shown in~\cite{JNS2} and~\cite{DJKNS1}. The
non-metacyclic maps ${\cal M}={\cal N}(2^e;k,l)$ in ${\cal R}(2^e)$ have mutually non-isomorphic
orientation-preserving automorphism groups ${\rm Aut}^+{\cal M}$, so they are invariant under $H_j$ for
each odd $j$. 

The operation $H_{-1}$ is of particular interest. An orientably regular map $\cal M$ is said to be {\sl
reflexible\/} if it has an automorphism reversing the orientation, or equivalently ${\cal M}$ is
isomorphic (as an oriented map) to its mirror image $\overline{\cal M}=H_{-1}({\cal M})$; otherwise
$\cal M$ and $\overline{\cal M}$ form a {\sl chiral pair}. If a map ${\cal M}\in{\cal R}(n)$ corresponds
to an isobicyclic triple $(G, x, y)\in{\cal I}(n)$, then $\overline{\cal M}$ corresponds to the triple
$(G, x^{-1}, y^{-1})$, so $\cal M$ is reflexible if and only if $G$ has an automorphism inverting $x$ and
$y$; we will then say that $(G, x, y)$ is reflexible. Similarly, if $\cal M$ corresponds to an
isobicyclic labelling $(\Gamma,\Lambda)\in{\cal L}(n)$, then $\overline{\cal M}$ corresponds to the
labelling $(\Gamma,\overline{\Lambda})$ where $\overline{\Lambda}(p)=(P, x_p^{-1}, y_p^{-1})$ if a vertex
$p$ has label $\Lambda(p)=(P, x_p, y_p)$, and $\overline{\Lambda}(q\rup p)=\Lambda(q\rup p)^{-1}$ for
each arc $q\rup p$. It follows that $\cal M$ is reflexible if and only if $\Lambda(p)$ is reflexible for
each vertex $p$ of $\Gamma$ and $\Lambda(q\rup p)^2=1$ for each arc $q\rup p$ of $\Gamma$. This allows
us to classify and to count the reflexible embeddings by considering their corresponding labellings.
 
First we consider vertex labels $\Lambda(p)$. It is shown in~\cite{JNS2} that for an odd prime $p$ the
only reflexible map in ${\cal R}(p^e)$ is the standard embedding, corresponding to an abelian Sylow
$p$-subgroup $P\cong C_{p^e}\times C_{p^e}$ of $G$, so the only reflexible label
$\Lambda(p)\in{\cal I}(p^e)$ is the standard triple $(P, x_p, y_p)$. In the case $p=2$, however,
there are extra possibilities. Firstly, as shown in~\cite{DJKNS1}, the reflexible metacyclic maps in ${\cal
R}(2^e)$ are the standard embedding ${\cal S}(2^e)={\cal M}(2^e,e,1)$, with $P\cong C_{2^e}\times
C_{2^e}$, and also (provided $e\geq 3$) the map ${\cal M}(2^e,e-1,1)$, corresponding to a nonabelian
group $P$ with $f=e-1$. In addition, the non-metacyclic maps ${\cal N}(2^e;k,l)$ in ${\cal R}(2^e)$ are
all reflexible: there is one such map if $e=2$, and there are four for each $e\geq 3$. Thus when $n$ is
even the number of possibilities for a reflexible label $\Lambda(2)$ is $1, 2$ or $6$ as $2^e\parallel n$
with $e=1$, $e=2$ or $e\geq 3$ respectively.

We now consider which arc labels $\lambda=\Lambda(q\rup p)$ can satisfy $\lambda^2=1$. By Lemma~7.4 we have $\lambda^{q^m}=1$ for some $m$, so if $q$ is odd then the only solution is $\lambda=1$; this corresponds to the trival action of a Sylow $q$-subgroup $Q$ on a Sylow $p$-subgroup $P$, so it cannot be the label of an arc in $\Gamma$. It follows that if $\cal M$ is reflexible then any arc in $\Gamma$ must have the form $2\rup p$. In particular, if $n$ is odd then $\Gamma$ must be the null graph; each vertex label $\Lambda(p)$ is a standard triple, since $p$ is odd, so in this case the standard embedding is the only reflexible map in ${\cal R}(n)$. If $n$ is even then any arc $2\rup p$ in $\Gamma$ must be labelled with $\lambda=-1$, since this is the only solution of $\lambda^2=1\neq \lambda$ when $p$ is odd, and conversely Lemma~7.4 shows that each of the reflexible $2$-subgroups $Q$ listed above can act diagonally on $P$ with eigenvalue $\lambda=-1$. It follows that if $n$ is divisible by $r$ odd primes $p$ then each of the $2^r$ possible subgraphs $\Gamma\preceq\Pi_n$ (one for each choice of a set of arcs $2\rup p$) has a unique isobicyclic arc-labelling. Taking account of the number of possible
vertex-labellings, considered above, we see that the number $\rho(n)$ of reflexible embeddings of
$K_{n,n}$ is given by
\[\rho(n)=\left\{
\begin{array}{ll}
1&\mbox{if $e=0$,}\\
2^r&\mbox{if $e=1$,}\\
2^{r+1}&\mbox{if $e=2$,}\\
3.2^{r+1}&\mbox{if $e\geq 3$,}
\end{array}
\right.
\]
where
\[n=2^ep_1^{e_1}\ldots p_r^{e_r}\]
for distinct odd primes $p_i$ and integers $e_i\geq 1$. This agrees with the enumeration recently
obtained by Kwak and Kwon in~\cite{KK2}, using a different method. The remaining maps in ${\cal R}(n)$ occur in chiral pairs, and the number $\chi(n)$ of such pairs can now be obtained from the formula
\[\chi(n)=\frac{1}{2}\bigl(\nu(n)-\rho(n)\bigr),\]
where $\nu(n)$ is the number $|{\cal R}(n)|$ of regular embeddings of $K_{n,n}$, to be determined in the
next two sections.

The above argument also gives an explicit description of the reflexible embeddings $\cal M$ of
$K_{n,n}$. As we have seen, if $n$ is odd then ${\cal M}$ is the standard embedding, so we may assume
that $n$ is even. In the notation of Corollary~4.2, some direct factor of $\cal S$, which we can
take to be ${\cal S}_1$, is one of the reflexible maps in ${\cal R}(2^e)$ described above, while
${\cal S}_2,\ldots,{\cal S}_k$, which are direct factors of $\cal M$, are all standard embeddings.
In the canonical decomposition, $S$ acts diagonally on $T$ with eigenvalue $-1$ (so that $S_2,\ldots,
S_k$, having odd orders, act trivially), and
\[{\cal M}\cong({\cal T}:{\cal S}_1)\times{\cal
S}_2\times\cdots\times{\cal S}_k\cong({\cal T}:{\cal S}_1)\times{\cal S}(n/2^et).\eqno(10.1)\]

The {\sl Petrie dual\/} ${\rm P}({\cal M})$ of a map $\cal M$, also discussed in~\cite{Wil}, has the same
embedded graph as $\cal M$, and has new faces, bounded by the Petrie polygons (closed zig-zag paths) of $\cal M$, so that ${\rm Aut}\,{\cal M}={\rm Aut}\,{\rm P}({\cal M})$. If ${\cal M}\in{\cal R}(n)$ then
${\rm P}({\cal M})$ is an embedding of $K_{n,n}$, and is orientable since this graph is bipartite, but it
need not be regular; if it is, then $\cal M$ must be reflexible, since a half-turn reversing an edge of
${\rm P}({\cal M})$ acts as a reflection reversing the same edge of $\cal M$. We will therefore consider
the effect of the Petrie duality operation $\rm P$ on the reflexible maps ${\cal M}\in{\cal R}(n)$.

If such a map $\cal M$ corresponds to a triple $(G, x, y)\in{\cal I}(n)$ then ${\rm P}({\cal M})$
corresponds to the triple $(G, x, y^{-1})$: this is in ${\cal I}(n)$ since the automorphism of $G$
transposing $x$ and $y$, composed with the automorphism inverting them, gives an automorphism
tranposing $x$ and $y^{-1}$. The Petrie polygons of $\cal M$ are $2m'$-gons, so ${\rm P}({\cal M})$ has
type $\{2m',n\}$, where $m'$ is the order of $xy^{-1}$. By arguments similar to those used in the case
of the operation $H_{-1}$, if $\cal M$ corresponds to a pair $(\Gamma,\Lambda)\in{\cal L}(n)$ then ${\rm
P}({\cal M})$ corresponds to the pair $(\Gamma,{\rm P}(\Lambda))$, where the labelling ${\rm
P}(\Lambda)$ is obtained from $\Lambda$ by replacing each vertex label $\Lambda(p)=(P, x_P, y_P)$ with
$(P, x_P, y^{-1}_P)$; the arc labels $\lambda=\Lambda(q\rup p)$ all satisfy $\lambda^{-1}=\lambda$, so
they are unchanged.
 
We say that a map $\cal M$ is {\sl self-Petrie\/} if ${\cal M}\cong{\rm P}({\cal M})$. A map ${\cal
M}\in{\cal R}(n)$ is self-Petrie if and only if $G$ has an automorphism fixing $x$ and inverting $y$, or
equivalently the vertex labels $\Lambda(p)$ are all self-Petrie. Since every standard embedding is
self-Petrie, it follows that if $\cal M$ is as in $(10.1)$ then
\[{\rm P}({\cal M})\cong({\cal T}:{\rm P}({\cal S}_1))\times{\cal
S}_2\times\cdots\times{\cal S}_k\cong({\cal T}:{\rm P}({\cal S}_1))\times{\cal S}(n/2^et).\]
Thus $\cal M$ is self-Petrie if and only if ${\cal S}_1$ is self-Petrie. Of the reflexible maps
in ${\cal R}(2^e)$ described above, all are self-Petrie with the exception of ${\cal N}(2^e\,;0,1)$ and
${\cal N}(2^e\,;1,1)$, which are Petrie duals of each other for each $e\geq 3$. (Note that ${\rm
Aut}\,{\cal N}(2^e\,;0,1)\cong{\rm Aut}\,{\cal N}(2^e\,;1,1)$, even though
${\rm Aut}^+{\cal N}(2^e\,;0,1)\not\cong{\rm Aut}^+{\cal N}(2^e\,;1,1)$.) It follows that a reflexible map
${\cal M}\in{\cal R}(n)$ is self-Petrie if and only if
${\cal S}_1$ is not one of these two maps, so the number $\sigma(n)$ of self-Petrie maps in ${\cal I}(n)$
is given by
\[\sigma(n)=\left\{\begin{array}{ll}
1&\mbox{if $e=0$,}\\
2^r&\mbox{if $e=1$,}\\
2^{r+1}&\mbox{if $e=2$,}\\
2^{r+2}&\mbox{if $e\geq 3$,}
\end{array}\right.\]
where $2^e\parallel n$ and $n$ is divisible by $r$ distinct odd primes. As in the case of
reflexibility, this agrees with the enumeration obtained by Kwak and Kwon in~\cite{KK2}.

\section{Enumeration: $n$ odd}

Having used the bijections between ${\cal R}(n)$, ${\cal I}(n)$ and ${\cal L}(n)$ to count the reflexible
and self-Petrie embeddings of $K_{n,n}$, we can use the same method to count all its regular embeddings.
By Theorem~2.1, the mapping ${\cal M}\mapsto(G, x, y)$ induces a bijection ${\cal R}(n)\to{\cal I}(n)$, so
the number
\[\nu(n)=|{\cal R}(n)|\]
of regular embeddings of $K_{n,n}$ is equal to the number $|{\cal I}(n)|$ of $n$-isobicyclic triples (up
to isomorphism in each case). By Proposition~8.1 the mapping $(G, x, y)\mapsto(\Gamma,\Lambda)$ is a
bijection ${\cal I}(n)={\cal J}(n)\to{\cal L}(n)$, so $\nu(n)$ is also equal to the number $|{\cal
L}(n)|$ of isobicyclic labellings of short spanning subgraphs $\Gamma\preceq\Pi_n$. We can therefore
enumerate maps by enumerating labellings. We will do this first in the simplest case, when $n$ is odd.

Given $\Gamma$, if $q\in N(\Gamma)$ and $q>2$, then by~\cite{JNS2} there are $\phi(q^{e_q-f})$ possible
labels $\Lambda(q)=(Q, x_q, y_q)$ for $q$ for each $f=1,\ldots, e_q$, where $q^{e_q}\parallel n$; given
any $p\in\Pi_n$ such that $q\rup p$, Lemma~7.4 implies that there are $q^{f_{q,p}}$ possible
diagonal actions of $Q$ on the abelian group $P$, where $f_{q,p}=\min\{f, r(q,p)\}$ and we define
$r=r(q,p)$ by $q^r\parallel p-1$. One of these is the trivial action, so there are $q^{f_{q,p}}-1$
nontrivial actions of $Q$ on $P$. If follows that there are
\[\sum_{f=1}^{e_q}\bigl(\phi(q^{e_q-f})\prod_{q\rup p}(q^{f_{q,p}}-1)\bigr)\]
possible labels for $q$ and its incident arcs $q\rup p$ in $\Gamma$. If $q\in T(\Gamma)$ then Lemma~7.5
implies that there is only one possibility for $Q$ and its actions, namely that
$Q$ is abelian and that it acts trivially on all $P$. It follows that if $n$ is odd then there are
\[\prod_{q\in N(\Gamma)}\Bigl(\sum_{f=1}^{e_q}\bigl(\phi(q^{e_q-f})
\prod_{q\rup p}(q^{f_{q,p}}-1)\bigr)\Bigr)\]
isobicyclic labellings $\Lambda$ of each short spanning subgraph $\Gamma\preceq\Pi_n$. Using the
bijection ${\cal R}(n)\to{\cal L}(n)$ described above we therefore have:

%Theorem 11.1.

\begin{thm}
 If $n$ is odd then
\[\nu(n)=\sum_{\Gamma\preceq\Pi_n}\Biggl(\prod_{q\in
N(\Gamma)}\Bigl(\sum_{f=1}^{e_q}\bigl(\phi(q^{e_q-f})
\prod_{q\rup p}(q^{f_{q,p}}-1)\bigr)\Bigr)\Biggr).\eqno{\square}\]
\end{thm}

For each $\Gamma\preceq\Pi_n$, the corresponding summand in this formula is nonnegative, so it provides a
lower bound for $\nu(n)$. Taking $\Gamma$ to be the null graph, so that $N(\Gamma)=\Pi_n$ and
$\prod_{q\rup p}(q^{f_{q,p}}-1)$ is the empty product, equal to $1$ for all primes $q$ dividing $n$, we
see that
\[\nu(n)\geq\prod_{q|n}\Bigl(\sum_{f=1}^{e_q}\phi(q^{e_q-f})\Bigr)=\prod_{q|n}q^{e_q-1}=n/\prod_{q|n}q.\]
This lower bound, first obtained in~\cite{JNS2}, represents the number of embeddings in ${\cal R}(n)$ for
which $G$ is nilpotent. It is attained if and only if $\Pi_n$ is itself a null graph, that is, no primes
$p$ and $q$ dividing $n$ satisfy $p\equiv 1$ mod~$(q)$. In particular, an odd integer $n$ satisfies
$\nu(n)=1$ if and only if $\Pi_n$ is null and $e_q=1$ for each prime $q$ dividing $n$, or equivalently
$\gcd(n,\phi(n))=1$, a result first obtained in~\cite{JNS1}.

The formula in Theorem~11.1 simplifies considerably when $n$ is divisible by a small number of primes. In~\cite{JNS2} it was shown that $\nu(p^e)=p^{e-1}$ for each odd prime $p$ (see also~\cite{NSZ} and~\cite{KK1} for the cases $e=1$ and $e=2$), so here we deal with the case where $n$ is divisible
by two odd primes. We considered the type and genus of the corresponding maps in Example~9.2.

% Corollary 11.2.

\begin{cor}
Let $n=p_1^dp_2^e$ where $p_1$ and $p_2$ are odd primes, $d, e\geq 1$ and
$p_1>p_2$. Then
$$\nu(n)=p_1^{d-1}p_2^{e-1}+\min\{e, r\}(p_2^e-p_2^{e-1})$$
where $p_2^r\parallel p_1-1$.
\end{cor}

\noindent{\sl Proof.} In this case $\Pi_n$ has two vertices $p_1$ and $p_2$, with an arc $p_2\to p_1$
if and only if $r\geq 1$, so each $\Gamma\preceq\Pi_n$ is either the null graph or equal to $\Pi_n$, the
latter possible only if $r\geq 1$.

If $\Gamma$ is the null graph then $q\in N(\Gamma)$ for both $q=p_1$ and $q=p_2$, and in each case
$\prod_{q\rup p}(q^{f_{q,p}}-1)$ is the empty product, equal to $1$, so the number of embeddings
corresponding to $\Gamma$ is
\[\sum_{f=1}^d\phi(p_1^{d-f})\sum_{f=1}^e\phi(p_2^{e-f})
=p_1^{d-1}p_2^{e-1}.\eqno(11.1)\]

If $r=0$, so that $p_2$ does not divide $p_1-1$, these are the only embeddings, but if $r\geq 1$
there are additional embeddings corresponding to $\Gamma=\Pi_n$. In this case $N(\Gamma)$ consists of a single prime $q=p_2$, with $q\rup p$ if and only if $p=p_1$, so the number of embeddings corresponding to $\Gamma$ is
\[\sum_{f=1}^e\phi(p_2^{e-f})(p_2^{f_{p_2,p_1}}-1)\eqno(11.2)\]
where $f_{p_2, p_1}=\min\{f,r\}$ and $p_2^r\parallel p_1-1$. If $r\geq e$ (equivalently $p_2^e\mid p_1-1$)
then $f_{p_2,p_1}=f$ for each $f=1,\ldots, e$, so the number of embeddings obtained in $(11.2)$ is
\begin{eqnarray*}
\sum_{f=1}^e\phi(p_2^{e-f})(p_2^f-1)& = &\sum_{f=1}^e\phi(p_2^{e-f})p_2^f-\sum_{f=1}^e\phi(p_2^{e-f})\\
& = &\sum_{f=1}^{e-1}(p_2^{e-f}-p_2^{e-f-1})p_2^f+p_2^e-p_2^{e-1}\\
& = &(e-1)(p_2^e-p_2^{e-1})+p_2^e-p_2^{e-1}\\
& = &e(p_2^e-p_2^{e-1}).
\end{eqnarray*}
%\eqno(11.3)
If $r<e$ then $f_{p_2,p_1}=f$ for each $f=1,\ldots, r$, and $f_{p_2,p_1}=r$ for each $f=r,\ldots,
e$; for each $f=r+1,\ldots, e$ we must therefore replace the summand $\phi(p_2^{e-f})(p_2^f-1)$
in the calculation $(11.3)$ with $\phi(p_2^{e-f})(p_2^r-1)$, thus subtracting
\begin{eqnarray*}
\sum_{f=r+1}^e\phi(p_2^{e-f})(p_2^f-p_2^r)
&=&\sum_{f=r+1}^{e-1}(p_2^{e-f}-p_2^{e-f-1})(p_2^f-p_2^r)+p_2^e-p_2^r\\
&=&(e-1-r)(p_2^e-p_2^{e-1})-\sum_{f=r+1}^{e-1}(p_2^{e-f+r}-p_2^{e-f-1+r})+p_2^e-p_2^r\\
&=&(e-1-r)(p_2^e-p_2^{e-1})-(p_2^{e-1}-p_2^r)+p_2^e-p_2^r\\
&=&(e-r)(p_2^e-p_2^{e-1})
\end{eqnarray*}

\noindent from the final total, so in this case the number of embeddings obtained in $(11.2)$ is
\[e(p_2^e-p_2^{e-1})-(e-r)(p_2^e-p_2^{e-1})=r(p_2^e-p_2^{e-1}).\eqno(11.4)\]
Combining the two possibilities $r\geq e$ and $r<e$ covered by $(11.3)$ and $(11.4)$, we see that the
number of embeddings corresponding to $\Gamma=\Pi_n$ in $(11.2)$ is
\[\min\{e, r\}(p_2^e-p_2^{e-1}).\]
Adding this to $(11.1)$ we find that the total number of embeddings is
\[\nu(n)=p_1^{d-1}p_2^{e-1}+\min\{e, r\}(p_2^e-p_2^{e-1}).\eqno{\square}\]

\medskip

The embeddings enumerated by the first summand $p_1^{d-1}p_2^{e-1}$, corresponding to the null graph $\Gamma$, are simply the cartesian products of the maps in ${\cal R}(p_1^d)$ and ${\cal
R}(p_2^e)$. These embeddings have $\pi=\emptyset$ in Proposition~9.2, so they have type $\{2n,n\}$
and genus
\[g=\frac{(n-1)(n-2)}{2}.\]
If $r\geq 1$ then the embeddings enumerated by the second summand have $\pi=\{p_1\}$, so by
Proposition~9.2 they have type $\{2p_2^e, n\}$ and genus
\[g=1+\frac{n}{2}\Bigl(n-p_1^d-2\Bigr).\]

Corollary~11.2, together with the results for odd prime powers obtained in~\cite{JNS2} (see \S3), covers all odd $n\leq 120$ with the exception of $n=105=3.5.7$, where Theorem~11.1 is required; the corresponding values of $\nu(n)$ are given in Table~1, together with those for even $n$ to be considered in \S12.
 
Kwak and Kwon~\cite{KK1} have developed a different method for studying the regular embeddings of
$K_{n,n}$, based on describing them in terms of permutations in the symmetric group $S_n$. For instance, they deal with the case where $n$ is a product of two primes, obtaining the enumeration in Corollary~11.2 in the case $d=e=1$. They have used this method to carry out a computer analysis of the maps in ${\cal R}(n)$ for all $n\leq 20$, and Fujisaki~\cite{Fuj} has extended their enumeration to $n\leq 100$. The values for $\nu(n)$ given here agree with his for all odd $n$ in this range.

\[
\begin{array}{ccccccccccccccccc}
n & \nu(n) && n & \nu(n) && n & \nu(n) && n & \nu(n) && n & \nu(n) && n & \nu(n)\\
1&1 && 21&3 && 41 & 1 && 61 & 1 && 81 & 27 && 101 &  $\;$1\\
2&1 && 22&2 && 42 & 8 && 62 & 2 && 82 & 2 && 102 & $\;$4\\
3 & 1 && 23 & 1 && 43 & 1 && 63 & 9 && 83 & 1 && 103 & $\;$1\\
4 & 2 && 24 & 12 && 44 & 4 && 64 & 20 && 84 & 16 && 104 &  $\;$16\\
5 & 1 && 25 & 5 && 45 & 3 && 65 & 1 && 85 & 1 && 105 & $\;$3\\
6 & 2 && 26 & 2 && 46 & 2 && 66 & 4 && 86 & 2 && 106 & $\;$2\\
7 & 1 && 27 & 9 && 47 & 1 && 67 & 1 && 87 & 1 && 107 & $\;$1\\
8 & 6 && 28 & 4 && 48 & 16 && 68 & 6 && 88 & 12 && 108 & $\;$20\\
9 & 3 && 29 & 1 && 49 & 7 && 69 & 1 && 89 & 1 && 109 & $\;$1\\
10 & 2 && 30 & 4 && 50 & 6 && 70 & 4 && 90 & 8 && 110 & $\;$12\\
11 & 1 && 31 & 1 && 51 & 1 && 71 & 1 && 91 & 1 && 111 & $\;$3\\
12 & 4 && 32 & 12 && 52 & 6 && 72 & 24 && 92 & 4 && 112 & $\;$16\\
13 & 1 && 33 & 1 && 53 & 1 && 73 & 1 && 93 & 3 && 113 & $\;$1\\
14 & 2 && 34 & 2 && 54 & 10 && 74 & 2 && 94 & 2 && 114 & $\;$8\\
15 & 1 && 35 & 1 && 55 & 5 && 75 & 5 && 95 & 1 && 115 & $\;$1\\
16 & 8 && 36 & 8 && 56 & 12 && 76 & 4 && 96 & 24 && 116 & $\;$6\\
17 & 1 && 37 & 1 && 57 & 3 && 77 & 1 && 97 & 1 && 117 & $\;$9\\
18 & 4 && 38 & 2 && 58 & 2 && 78 & 8 && 98 & 8 && 118 & $\;$2\\
19 & 1 && 39 & 3 && 59 & 1 && 79 & 1 && 99 & 3 && 119 & $\;$1\\
20 & 6 && 40 & 16 && 60 & 12 && 80 & 24 && 100 & 14 && 120 & $\;$32
\end{array}
\]
\smallskip

\centerline{Table 1: values of $\nu(n)=|{\cal R}(n)|$ for $1\leq n\leq 120$}

\section{Enumeration: $n$ even}

The method of enumeration used in \S11 can also be used when $n$ is even, so that $2\in N(\Gamma)$ for all
$\Gamma\preceq\Pi_n$. Suppose that $2^{e_2}\parallel n$, with $e_2\geq 1$. The following result shows
that the formula in Theorem~11.1 for $n$ odd is in fact valid for all $n$ not divisible by $8$.

% Theorem 12.1.

\begin{thm}
If $1\leq e_2\leq 2$, so that $2\parallel n$ or $2^2\parallel n$, then
\[\nu(n)=\sum_{\Gamma\preceq\Pi_n}\Biggl(\prod_{q\in
N(\Gamma)}\Bigl(\sum_{f=1}^{e_q}\bigl(\phi(q^{e_q-f})
\prod_{q\rup p}(q^{f_{q,p}}-1)\bigr)\Bigr)\Biggr).\]
\end{thm}

\medskip

\noindent{\sl Proof.} If $e_2=1$ then the label $\Lambda(2)$ attached to the vertex $q=2$ of $\Pi_n$ must
be the standard triple in ${\cal I}(2)$, and the Sylow $2$-subgroup $Q\cong C_2\times C_2$ of $G$ has two
possible actions on an abelian Sylow $p$-subgroup $P$ for $p>2$, with eigenvalues $\pm 1$.
This is consistent with the formula in Theorem~11.1, with
$f=1$ and the number of actions given by $2^{f_{2,p}}$. 

If $e_2=2$ there is one metacyclic Sylow $2$-subgroup $Q$ which can occur, namely $C_4\times C_4$,
corresponding to $f=2$; unlike in the odd prime power case, there is no metacyclic group for $f=1$,
but instead there is one non-metacyclic group $Q=G(4;0,0)$ which can occur. The number of possible
actions of $C_4\times C_4$ on $P$ is $2^{f_{2,p}}$, and for $G(4;0,0)$ it is $2$, which is the value of
$2^{f_{2,p}}$ when $f=1$. Thus the formula in Theorem~11.1 remains valid when $e_2=2$.
\hfill$\square$

\medskip

Similar arguments show that Corollary~11.2 also applies when $p_2=2$ and $e=1$ or $2$. If $e=1$, so that
$n=2p_1^d$ for an odd prime $p_1$, it gives $\nu(n)=p_1^{d-1}+1$, and if $e=2$, so that $n=4p_1^d$, it
gives $\nu(n)=2p_1^{d-1}+4$ or $2p_1^{d-1}+2$ as $p_1\equiv 1$ or $-1$ mod~$(4)$.

The situation is more complicated when $n$ is divisible by a higher power of $2$:

% Theorem 12.2.

\begin{thm}
If $n$ is divisible by $8$ then
\[\nu(n)=\negthinspace\sum_{\Gamma\preceq\Pi_n}\Biggl(\Bigl(\sum_{f=2}^{e_2}\bigl(\phi(2^{e_2-f})
\prod_{2\rup p}(2^{f_{2,p}}-1)\bigr)+4\Bigr)
\prod_{2\neq q\in N(\Gamma)}\Bigl(\sum_{f=1}^{e_q}\bigl(\phi(q^{e_q-f})
\prod_{q\rup p}(q^{f_{q,p}}-1)\bigr)\Bigr)\Biggr).\]

% Alternative broken display:
%\begin{eqnarray*}
%\nu(n)&=&\negthinspace\sum_{\Gamma\preceq\Pi_n}\Biggl(\Bigl(\sum_{f=2}^{e_2}\bigl(\phi(2^{e_2-f})
%\prod_{2\rup p}(2^{f_{2,p}}-1)\bigr)+4\Bigr)\\
%&&\qquad\qquad\times\prod_{2\neq q\in N(\Gamma)}\Bigl(\sum_{f=1}^{e_q}\bigl(\phi(q^{e_q-f})
%\prod_{q\rup p}(q^{f_{q,p}}-1)\bigr)\Bigr)\Biggr).
%\end{eqnarray*}
\end{thm}

An equivalent version of this formula is:
\begin{eqnarray*}
\nu(n)&=&\negthinspace\sum_{\Gamma\preceq\Pi_n}\Biggl(\prod_{q\in N(\Gamma)}
\Bigl(\sum_{f=1}^{e_q}\bigl(\phi(q^{e_q-f})\prod_{q\rup p}(q^{f_{q,p}}-1)\bigr)\Bigr)\Biggr)\\
&&+\,(4-2^{e_2-2})\sum_{\Gamma\preceq\Pi_n}\Biggl(\prod_{2\neq q\in N(\Gamma)}
\Bigl(\sum_{f=1}^{e_q}\bigl(\phi(q^{e_q-f})\prod_{q\rup p}(q^{f_{q,p}}-1)\bigr)\Bigr)\Biggr).
\end{eqnarray*}

\noindent{\sl Proof.} If $e_2\geq 3$ there are $2^{e_2-2}$ metacyclic $2$-groups $Q$ which can occur,
namely $\phi(2^{e_2-f})$ for each $f=2,\ldots, e_2$ (again, $f=1$ is excluded), and in addition there
are four non-metacyclic groups. In the metacyclic case, the number of possible actions of $Q$ on an
abelian Sylow $p$-subgroup $P$ is again $2^{f_{2,p}}$. In the non-metacyclic case, there are two
actions for each $Q$, with eigenvalues $\lambda=\pm 1$. It follows that, for a given
$\Gamma\preceq\Pi_n$, the number of possibilities for $Q$ and its (nontrivial) actions on the Sylow
$p$-subgroups such that $2\rup p$ is
\[\sum_{f=2}^{e_2}\bigl(\phi(2^{e_2-f})\prod_{2\rup p}(2^{f_{2,p}}-1)\bigr)+4.\]
When $f=1$ we have $f_{2,p}=1$ for all odd primes $p$, so the `missing' summand for $f=1$ is
\[\phi(2^{e_2-f})\prod_{2\rup p}(2^{f_{2,p}}-1)=\phi(2^{e_2-1})=2^{e_2-2}\]
and hence we can write the number of possibilities for $Q$ and its actions as
\[\sum_{f=1}^{e_2}\bigl(\phi(2^{e_2-f})\prod_{2\rup p}(2^{f_{2,p}}-1)\bigr)-2^{e_2-2}+4.\]
The main sum here is given by putting $q=2$ in the corresponding formula for odd primes $q$, and we
can regard $-2^{e_2-2}+4$ as a `correction term' for $e_2\geq 3$. \hfill$\square$

\medskip

This correction term is $0$ if $e_2=4$, so the formul\ae\/ in Theorem~11.1 and Corollary~11.2 are
also valid in this case. Thus if $n=16p_1^d$ for an odd prime $p_1$ then Corollary~11.2 implies that
$\nu(n)=8p_1^{d-1}+8\min\{4,r\}$ where $2^r\parallel p_1-1$. More generally, the analogue of
Corollary~11.2 for even $n$ is as follows:

\medskip

\noindent{\bf Example 12.3.} Let $n=2^ep_1^d$ where $p_1$ is an odd prime, $e\geq 3$, and $d\geq 1$. As
in Corollary~11.2, the first term takes the value
\[2^{e-1}p_1^{d-1}+\min\{e, r\}(2^e-2^{e-1})=2^{e-1}(p_1^{d-1}+\min\{e, r\})\]
where $2^r\parallel p_1-1$. In the second summation, if $\Gamma=\Pi_n$ then there is no odd $q$ in
$N(\Gamma)$, so the product is the empty product, with value $1$. If $\Gamma$ is the null subgraph
then $N(\Gamma)=\{p_1\}$, and there is no $p$ such that $p_1\rup p$, so
\begin{eqnarray*}
\nu(n)&=&2^{e-1}(p_1^{d-1}+\min\{e, r\})+(4-2^{e-2})(1+p_1^{d-1})\\
&=&(2^{e-2}+4)p_1^{d-1}+(2\min\{e,r\}-1)2^{e-2}+4.
\end{eqnarray*}
For instance, if $e=3$, so that $n=8p_1^d$, then
\[\nu(n)=\left\{
\begin{array}{ll}
6p_1^{d-1}+14&\mbox{if $p_1\equiv 1$ mod~$(8)$,}\\
&\\
6p_1^{d-1}+10&\mbox{if $p_1\equiv 5$ mod~$(8)$,}\\
&\\
6p_1^{d-1}+6&\mbox{if $p_1\equiv 3$ or $7$ mod~$(8)$.}
\end{array}
\right.
\]
Similarly, if $e=4$, so that $n=16p_1^d$, then
\[\nu(n)=\left\{
\begin{array}{ll}
8p_1^{d-1}+32&\mbox{if $p_1\equiv 1$ mod~$(16)$,}\\
&\\
8p_1^{d-1}+16&\mbox{if $p_1\equiv 5$ or $13$ mod~$(16)$,}\\
&\\
8p_1^{d-1}+8&\mbox{if $p_1\equiv 3, 7, 11$ or $15$ mod~$(16)$.}
\end{array}
\right.
\]
These results, together with those for $n=2^e$ in~\cite{DJKNS1, DJKNS2}, account for most of the
entries for even $n\leq 120$ in Table~1; the remaining entries can be obtained from Theorem~12.1. For
$n\leq 20$ the resulting values of $\nu(n)$ all agree with those in~\cite{KK1}, and for $n\leq 100$ with
those in~\cite{Fuj}, except that the latter gives $\nu(90)=6$. This may simply be a transcription error, as
there are seven rather obvious embeddings for $n=90$ which decompose as cartesian products, together with one indecomposable embedding. This is the only integer $n\leq 100$ for which our value for $\nu(n)$
differs from that in~\cite{Fuj}.

\section{The directed graphs $\Pi_n$}

The vertices of the directed graph $\Pi_n$ are the primes $p$ dividing $n$, with an arc $q\to p$
whenever $q$ divides $p-1$. We can label each such arc with the positive integer $r=r(q,p)$ such that
$q^r\parallel p-1$. In view of the preceding enumerative results it is interesting to know which labelled
directed graphs can occur in this context. A directed graph is said to be {\sl acyclic\/} if it contains
no directed cycles. An {\sl isomorphism\/} of labelled directed graphs is an isomorphism of directed
graphs such that each arc has the same label as its image.

% Proposition 13.1.

\begin{prop}
Let $\Delta$ be a finite directed graph whose arcs are labelled with positive
integers. Then $\Delta$ is isomorphic as a labelled directed graph to $\Pi_n$ for some $n$ if and only if
$\Delta$ is acyclic.
\end{prop}

\noindent{\sl Proof.} Each $\Pi_n$ is acyclic, since if it has an arc $q\to p$ then $q<p$; hence any directed graph isomorphic to $\Pi_n$ is also acyclic.

For the converse, we use induction on the number $l$ of vertices of $\Delta$. For technical reasons, we
will prove that if $\Delta$ is acyclic then $\Delta\cong\Pi_n$ for some odd $n$. This result is
trivial if $l=1$, so suppose that $l>1$ and the result has been proved for all labelled directed graphs
with $l-1$ vertices. Being finite and acyclic, $\Delta$ has a vertex $v$ which is not the source of
any arc (for instance, the last vertex in a maximal directed path). Let $\Delta'$ be the labelled
directed graph formed from $\Delta$ by removing $v$ and any incident arcs $u\to v$. Then $\Delta'$ is
acyclic, so the induction hypothesis implies that there is an isomorphism of labelled directed graphs
$\theta:\Delta'\to\Pi_{n'}$ for some odd integer $n'$. Each vertex $w$ of $\Delta'$ corresponds to a prime
$\theta(w)=q>2$ dividing $n'$. We will extend $\theta$ to an isomorphism $\Delta\to\Pi_n$ by taking
$n=n'p$ where $p$ is a suitably chosen odd prime, and defining $\theta(v)=p$. Such an extension is an
isomorphism if and only if $q^{r_q}\parallel p-1$ for each prime $q$ dividing $n'$, where $r_q=r$ if there
is an arc $w\to v$ in $\Delta$ labelled $r$, and $r_q=0$ if there is no arc $w\to v$. These conditions are
all satisfied if $p\equiv q^{r_q}+1$ mod~$(q^{r_q+1})$ for each $q$. The Chinese Remainder Theorem implies that this set of congruences is equivalent to a single congruence of the form $p\equiv a$ mod~$(b)$, where $b$ is the product of the prime powers $q^{r_q+1}$. Since $q^{r_q}+1$ is coprime to $q^{r_q+1}$ for each $q$ it follows that $a$ and $b$ are coprime, so Dirichlet's Theorem on primes in arithmetic progressions implies that there is at least one odd prime $p$ satisfying this congruence (in fact, there are
infinitely many).
\hfill$\square$

\medskip

A similar argument shows that the directed graphs corresponding to even integers $n$ are those which are
acyclic and have a vertex $u$ (corresponding to the prime $2$) with an arc from $u$ to every other vertex.

% Corollary 13.2.

\begin{cor}
Every finite graph is isomorphic to the underlying graph of some directed
graph $\Pi_n$.
\end{cor}

\noindent{\sl Proof.} If we number the vertices of a finite graph $v_1, v_2,\ldots$, and regard any edge
$v_iv_j$ as an arc $v_i\to v_j$ where $i<j$, we obtain an acyclic directed graph, so the result follows
from Proposition~13.1.
\hfill$\square$

\section{Connections with the random graph}

The arguments used in the preceding section can be adapted to establish a link between the directed
graph $\Pi$ and the {\sl random graph\/} or {\sl universal graph\/} $R$ studied by Erd\H os and
R\'enyi~\cite{ER} and constructed by Rado~\cite{Rad}. This graph has many remarkable properties, described in some detail in~\cite[\S5.1]{Cam} and~\cite[\S9.6]{DM}. For instance, given a countably infinite vertex set, if pairs of vertices are chosen randomly to be edges, each with probability~\(\frac{1}{2}\,\), then the resulting graph is isomorphic to $R$ with probability $1$.
 
In $\Pi$ there is an arc from the vertex $2$ to every other vertex. If we delete this vertex and all such
incident arcs, and ignore the direction of each remaining arc, we obtain an undirected graph $\Pi'$ whose
vertices are the odd primes, with an edge between $p$ and $q$ if and only if $q$ divides $p-1$ or vice
versa.

% Proposition 14.1.

\begin{prop}
$\Pi'$ is isomorphic to $R$.
\end{prop}

\noindent{\sl Proof.} The graph $R$ is characterised (up to isomorphism) among countable graphs by the
property that for each disjoint pair $U$ and $V$ of finite sets of vertices, there is a vertex adjacent to
every vertex in $U$ and to no vertex in $V$. The fact that $\Pi'$ has this property follows immediately
from the Chinese Remainder Theorem and Dirichlet's Theorem on primes in arithmetic progressions: the
simultaneous congruences $p\equiv 1$ mod~$(u)$ for all $u\in U$ and $p\equiv -1$ mod~$(v)$ for all $v\in V$ are equivalent to a single congruence modulo $\prod u\,.\prod v$, and this is satisfied by at least
one odd prime $p$ (by infinitely many, in fact). Then $p-1$ is divisible by each $u\in U$ and by no
$v\in V$, and $p$ cannot divide any $v-1$ (for otherwise $v=p+1$ would be even), so $\Pi'\cong R$.
\hfill $\square$

\medskip

Since every countable graph is isomorphic to an induced subgraph of $R$, this result gives an alternative
proof of Corollary~13.2.

\bigskip

\noindent{\bf Acknowledgements.} The author is grateful to Roman Nedela and Martin \v Skoviera for
their advice and their continued encouragement to undertake this research, and to J\"urgen Wolfart and the Mathematics Faculty at the J.~W.~Goethe University, Frankfurt-am-Main, for providing the ideal environment in which it could be carried out.

\end{document}